\documentclass[11pt,a4paper]{amsart} 
\usepackage[all]{xy}
\SelectTips{cm}{}
\usepackage{ifthen} 
\usepackage{amssymb}
\calclayout
\makeatletter 
\def\serieslogo@{} 
\def\@setcopyright{} 
\makeatother

\title{The stable derived category of a noetherian scheme}

\author{Henning Krause}
\address{Henning Krause\\ Institut f\"ur Mathematik\\
Universit\"at Paderborn\\ 33095 Paderborn\\ Germany.}
\email{hkrause@math.uni-paderborn.de}

\dedicatory{Dedicated to Claus Michael Ringel on the
occasion of his sixtieth birthday.}  


\newtheorem{lem}{Lemma}[section]
\newtheorem{prop}[lem]{Proposition} \newtheorem{cor}[lem]{Corollary}
\newtheorem{thm}[lem]{Theorem}

\theoremstyle{remark}

\theoremstyle{definition}
\newtheorem{exm}[lem]{Example}
\newtheorem{defn}[lem]{Definition}

\newtheorem{rem}[lem]{Remark}
\numberwithin{equation}{section}

\renewcommand{\mod}{\operatorname{mod}\nolimits}

\newcommand{\Flat}{\operatorname{Flat}\nolimits}
\newcommand{\proj}{\operatorname{proj}\nolimits}

\newcommand{\id}{\operatorname{id}\nolimits}

\newcommand{\Mod}{\operatorname{Mod}\nolimits}

\newcommand{\Hom}{\operatorname{Hom}\nolimits}

\newcommand{\END}{\operatorname{\mathcal E\!\!\:\mathit n\mathit d}\nolimits}
\newcommand{\HOM}{\operatorname{\mathcal H\!\!\:\mathit o\mathit m}\nolimits}
\renewcommand{\Im}{\operatorname{Im}\nolimits}
\newcommand{\Ker}{\operatorname{Ker}\nolimits}
\newcommand{\Coker}{\operatorname{Coker}\nolimits}

\newcommand{\Ext}{\operatorname{Ext}\nolimits}
\newcommand{\uExt}{\operatorname{\underline{Ext}}\nolimits}
\newcommand{\tExt}{\operatorname{\widehat{Ext}}\nolimits}

\newcommand{\noeth}{\operatorname{noeth}\nolimits}

\newcommand{\umod}{\operatorname{\underline{mod}}\nolimits}
\newcommand{\uMod}{\operatorname{\underline{Mod}}\nolimits}

\newcommand{\uHom}{\operatorname{\underline{Hom}}\nolimits}

\newcommand{\Supp}{\operatorname{Supp}\nolimits}

\newcommand{\Inj}{\operatorname{Inj}\nolimits}
\newcommand{\GInj}{\operatorname{GInj}\nolimits}
\newcommand{\Proj}{\operatorname{Proj}\nolimits}
\newcommand{\Qcoh}{\operatorname{Qcoh}\nolimits}
\newcommand{\coh}{\operatorname{coh}\nolimits}

\newcommand{\ac}{\mathrm{ac}}
\newcommand{\tac}{\mathrm{tac}}
\newcommand{\op}{\mathrm{op}}
\newcommand{\inc}{\mathrm{inc}}
\newcommand{\perf}{\mathrm{perf}}
\newcommand{\inj}{\mathrm{inj}}
\newcommand{\can}{\mathrm{can}}
\newcommand{\Id}{\mathrm{Id}}

\newcommand{\dg}{\mathrm{dg}}
\newcommand{\comp}{\mathop{\raisebox{+.3ex}{\hbox{$\scriptstyle\circ$}}}}
\newcommand{\lto}{\longrightarrow}
\newcommand{\xto}{\xrightarrow}

\def\a{\alpha}
\def\b{\beta}
\def\e{\varepsilon}
\def\d{\delta}

\def\p{\phi}
\def\r{\rho}
\def\s{\sigma}

\def\la{\lambda}

\def\Ga{\Gamma}
\def\La{\Lambda}
\def\Si{\Sigma}

\def\A{{\mathcal A}}
\def\B{{\mathcal B}}
\def\C{{\mathcal C}}

\def\K{{\mathcal K}}

\def\OO{{\mathcal O}}

\def\S{{\mathcal S}}
\def\X{{\mathcal X}}
\def\Y{{\mathcal Y}}

\def\T{{\mathcal T}}
\def\U{{\mathcal U}}

\def\bbP{\mathbb P}
\def\bbX{\mathbb X}
\def\bbY{\mathbb Y}
\def\bbZ{\mathbb Z}

\def\bfC{\mathbf C}
\def\bfD{\mathbf D}
\def\bfK{\mathbf K}
\def\bfL{\mathbf L}
\def\bfR{\mathbf R}
\def\bfS{\mathbf S}
\def\bfT{\mathbf T}
\begin{document}

\begin{abstract}
For a noetherian scheme, we introduce its unbounded stable derived
category. This leads to a recollement which reflects the passage from
the bounded derived category of coherent sheaves to the quotient
modulo the subcategory of perfect complexes. Some applications are
included, for instance an analogue of maximal Cohen-Macaulay
approximations, a construction of Tate cohomology, and an extension of
the classical Grothendieck duality. In addition, the relevance of the
stable derived category in modular representation theory is indicated.
\end{abstract}
\maketitle

\section{Introduction}

Let $\bbX$ be a separated noetherian scheme and denote by $\Qcoh\bbX$ the
category of quasi-coherent sheaves on $\bbX$. We consider the derived
category $\bfD(\Qcoh\bbX)$ and two full subcategories 
$$\bfD^\perf(\coh\bbX)\subseteq\bfD^b(\coh\bbX)\subseteq\bfD(\Qcoh\bbX)$$
which are of particular interest. Here, $\bfD^b(\coh\bbX)$ denotes the
bounded derived category of coherent sheaves, and $\bfD^\perf(\coh\bbX)$
denotes the subcategory of perfect complexes. 

Now let $\Inj\bbX$ be the full subcategory of injective objects in
$\Qcoh\bbX$, and denote by $\bfK(\Inj\bbX)$ its homotopy category. 
The composite
$$Q\colon\xymatrix{\bfK(\Inj\bbX)\ar[r]^-{\mathrm{inc}}&\bfK(\Qcoh\bbX)
\ar[r]^-{\mathrm{can}}&\bfD(\Qcoh\bbX)}$$ gives rise to a localization
sequence $$\xymatrix{\bfS(\Qcoh\bbX)\ar[r]^-I&\bfK(\Inj\bbX)
\ar[r]^-Q&\bfD(\Qcoh\bbX)}$$
where $\bfS(\Qcoh\bbX)$ denotes the full subcategory of all acyclic
complexes in $\bfK(\Inj\bbX)$. Thus $Q$ induces an equivalence
$$\bfK(\Inj\bbX)/\bfS(\Qcoh\bbX)\stackrel{\sim}\lto\bfD(\Qcoh\bbX).$$

Next we recall that an object $X$ in some category with coproducts is
{\em compact} if every map $X\to\coprod_iY_i$ into an arbitrary
coproduct factors through a finite coproduct.  For instance, an object
in $\bfD(\Qcoh\bbX)$ is compact if and only if it is isomorphic to a
perfect complex.  It is well known that the derived category
$\bfD(\Qcoh\bbX)$ is {\em compactly generated}, that is, there is a
set of compact objects which generate $\bfD(\Qcoh\bbX)$ \cite{N1}.  To
formulate our main result, let us denote by $\bfK^c(\Inj\bbX)$ and
$\bfS^c(\Qcoh\bbX)$ the full subcategories of compact objects in
$\bfK(\Inj\bbX)$ and $\bfS(\Qcoh\bbX)$ respectively.

\begin{thm} 
Let $\bbX$ be a separated noetherian scheme.
\begin{enumerate}
\item The functors $I,Q$ have left adjoints $I_\la,Q_\la$ and right
adjoints $I_\r,Q_\r$ respectively. We have therefore a
recollement
$$\xymatrix{\bfS(\Qcoh\bbX)\ar[r]&\bfK(\Inj\bbX)
\ar[r]\ar@<.75ex>[l]\ar@<-.75ex>[l]&\bfD(\Qcoh\bbX).
\ar@<.75ex>[l]\ar@<-.75ex>[l]}$$

\item The triangulated category $\bfK(\Inj\bbX)$ is compactly
generated, and $Q$ induces an equivalence
$\bfK^c(\Inj\bbX)\stackrel{}\to\bfD^b(\coh\bbX)$.
\item The sequence $$\xymatrix{\bfD(\Qcoh\bbX)\ar[r]^-{Q_\la}&
\bfK(\Inj\bbX)\ar[r]^-{I_\la}&\bfS(\Qcoh\bbX)}$$ is a localization
sequence. Therefore $\bfS(\Qcoh\bbX)$ is compactly generated, and
$I_\la\comp Q_\r$ induces (up to direct factors) an equivalence
$$\bfD^b(\coh\bbX)/\bfD^\perf(\coh\bbX)\stackrel{\sim}\lto\bfS^c(\Qcoh\bbX).$$
\end{enumerate}
\end{thm}

Note that this theorem is a special case of a general result about
Grothendieck categories. All we need is a locally noetherian
Grothendieck category $\A$, for instance $\A=\Qcoh\bbX$, such that
$\bfD(\A)$ is compactly generated. There is a surprising
consequence which seems worth mentioning.

\begin{cor} 
Let $\bbX$ be a separated noetherian scheme. Then a product of acyclic
complexes of injective objects in $\Qcoh\bbX$ is acyclic.
\end{cor}

We call the category $\bfS(\Qcoh\bbX)$ the {\em stable derived
category} of $\Qcoh\bbX$.  A first systematic study of the bounded
stable derived category
$$\bfD^b(\coh\bbX)/\bfD^\perf(\coh\bbX)$$ can be found in work of
Buchweitz \cite{Bu}. Unfortunately this beautiful paper has never been
published; see however \cite{BEH}.  He identifies for a Gorenstein
ring $\La$ the bounded derived category of finitely generated
$\La$-modules modulo perfect complexes
$$\bfD^b(\mod\La)/{\bfD^\perf(\mod\La)}$$ with the stable category of
maximal Cohen-Macaulay $\La$-modules and with the category of acyclic
complexes of finitely generated projective $\La$-modules. The same
identification appears in \cite{R} for selfinjective algebras and
plays an important role in modular representation theory of finite
groups; see also \cite{KV}.  The approach in the present paper differs
from that of Buchweitz substantially because we work in the unbounded
setting and we use injective objects instead of projectives. This has
some advantages.  For instance, we have in any Grothendieck category
enough injectives but often not enough projectives.  On the other
hand, we obtain a recollement in the unbounded setting which does not
exists in the bounded setting. In fact, the celebrated theory of
maximal Cohen-Macaulay approximations \cite{AB} is described as
`decomposition' \cite{AB} or `glueing' \cite{Bu}, but finds a natural
interpretation as `recollement' in the sense of \cite{BBD} if one
passes to the unbounded setting.  To be precise, the recollement
$$\xymatrix{\bfS(\Mod\La)\ar[r]&\bfK(\Inj\La)
\ar[r]\ar@<.75ex>[l]^-{I_\la}\ar@<-.75ex>[l]_-{I_\r}&\bfD(\Mod\La)
\ar@<.75ex>[l]^-{Q_\la}\ar@<-.75ex>[l]_-{Q_\r}}$$ induces for any
Gorenstein ring $\La$ the Gorenstein injective approximation functor 
$$T\colon\Mod\La\xto{\can}\bfD(\Mod\La)\xto{I_\la\comp
Q_\r}\bfS(\Mod\La)\xto{Z^0}\uMod\La$$ where $\uMod\La$ denotes the
stable category modulo injective objects. For any $\La$-module $A$,
the Gorenstein injective approximation $A\to TA$ is the `dual' of the
maximal Cohen-Macaulay approximation which is based on projective
resolutions. Let us stress again that this approach generalizes to any
locally noetherian Grothendieck category $\A$ provided that $\bfD(\A)$
is compactly generated. 

Next we explain the connection between Gorenstein injective
approximations and Tate cohomology. We fix a locally noetherian
Grothendieck category $\A$ and pass from the stable derived category
$\bfS(\A)$ to the full subcategory $\bfT(\A)$ of {\em totally acyclic}
complexes. An object in $\A$ is by definition {\em Gorenstein
injective} if it is of the form $\Ker (X^0\to X^1)$ for some $X$ in
$\bfT(\A)$.  The inclusion $G\colon\bfT(\A)\to\bfK(\Inj\A)$ has a left
adjoint $G_\la$. Given an object $A$ in $\A$ with injective resolution
$iA$, we may think of $G_\la iA$ as a complete injective resolution of
$A$.  This leads to the following definition of {\em Tate cohomology
groups}
$$\tExt^n_\A(A,B)=H^n\Hom_\A(A,G_\la iB)$$ for any
$A,B$ in $\A$ and $n\in\bbZ$. This cohomology theory is symmetric in
the sense that for any $A$ in $\A$, we have $\tExt^*_\A(A,-)=0$ iff
$\tExt^*_\A(-,A)=0$ iff $\tExt^0_\A(A,A)=0$. Let $\X$ denote the class
of all objects $A$ such that $\tExt^*_\A(A,-)$ vanishes, and let $\Y$
be the class of Gorenstein injective objects in $\A$.

\begin{thm}\label{th:tate}
  Let $\A$ be a locally noetherian Grothendieck category and suppose
  that $\bfD(\A)$ is compactly generated. 
\begin{enumerate}
\item $\X=\{A\in\A\mid\Ext^1_\A(A,B)=0$\textrm{ for all }$B\in\Y\}$.
\item $\Y=\{B\in\A\mid\Ext^1_\A(A,B)=0$\textrm{ for all }$A\in\X\}$.
\item Every object $A$ in $\A$ fits into exact sequences
\begin{equation*}
0\to Y_A\to X_A\to A\to 0\quad\textrm{and}\quad0\to A\to Y^A\to X^A\to
0
\end{equation*} 
in $\A$ with $X_A,X^A$ in $\X$ and $Y_A,Y^A$ in $\Y$.
\item $\X\cap\Y=\Inj\A$.
\end{enumerate}
\end{thm} 

After explaining some historical backround, let us mention more
recent work on stable derived categories. For instance, Beligiannis
develops a general theory of `stabilization' in the framework of
relative homological algebra \cite{B}, and J{\o}rgensen studies the
category of `spectra' for a module category \cite{J}. Also, Orlov
discusses the category
$$\bfD^b(\coh\bbX)/\bfD^\perf(\coh\bbX)$$
under the name
`triangulated category of singularities' and points out some
connection with the Homological Mirror Symmetry Conjecture \cite{O}.
In any case, our notation $\bfS(\Qcoh\bbX)$ reflects this terminology.

Our main results suggests that the homotopy category $\bfK(\Inj\bbX)$
deserves some more attention. We may think of this category as the
`compactly generated completion' of the category $\bfD^b(\coh\bbX)$.
In fact, the category $\coh\bbX$ of coherent sheaves carries a natural
DG structure and its derived category $\bfD_\dg(\coh\bbX)$ is
equivalent to $\bfK(\Inj\bbX)$. This follows from Keller's work
\cite{Ke} and complements a recent result of Bondal and van den Bergh
\cite{BB} which says that $\bfD(\Qcoh\bbX)$ is equivalent to
$\bfD_\dg(A)$ for some DG algebra $A$.

As another application of our main result, let us mention that the
adjoint pair of functors $\bfR f_*$ and $f^!$ which establish the
Grothendieck duality for a morphism $f\colon \bbX\to\bbY$ between
schemes \cite{Ha,N1}, can be extended to a pair of adjoint functors between
$\bfK(\Inj\bbX)$ and $\bfK(\Inj\bbY)$. 

\begin{thm}\label{th:duality}
Let $f\colon \bbX\to\bbY$ be a morphism between separated
noetherian schemes. Denote by $\bfR
f_*\colon\bfD(\Qcoh\bbX)\to\bfD(\Qcoh\bbY)$ the right derived direct
image functor and by $f^!$ its right adjoint. Then there is an adjoint
pair of functors $\hat\bfR f_*$ and $\widehat{f^!}$ making the
following diagram commutative.
$$\xymatrix{\bfD(\Qcoh\bbX)\ar[d]^{Q_\la}\ar[rr]^-{\bfR f_*}&&\bfD(\Qcoh\bbY)&
\bfD(\Qcoh\bbY)\ar[d]^{Q_\r}\ar[rr]^-{f^!}&&\bfD(\Qcoh\bbX)
\\
\bfK(\Inj\bbX)\ar[rr]^-{\hat\bfR f_*}&&\bfK(\Inj\bbY)\ar[u]_Q&
\bfK(\Inj\bbY)\ar[rr]^-{\widehat{f^!}}&&\bfK(\Inj\bbX)\ar[u]_Q\\
}$$
\end{thm}

Again, this theorem is really a lot more general. It is irrelevant
that the functor $f_*$ comes from a morphism $f\colon \bbX\to\bbY$.
All we need is that $f_*$ and its right derived functor $\bfR f_*$
preserve coproducts. On the other hand, there is a strengthened
version of Theorem~\ref{th:duality} which uses the special properties of
$f_*$. I am grateful to Amnon Neeman for pointing out the following.

\begin{thm}[Neeman]\label{th:sduality}
Let $f\colon \bbX\to\bbY$ be a morphism between separated
noetherian schemes. Then $\hat\bfR f_*$ sends acyclic complexes
to acyclic complexes. Thus we have an adjoint pair of functors
between $\bfS(\Qcoh\bbX)$ and $\bfS(\Qcoh\bbY)$, making
the following diagrams commutative.
$$\xymatrix{\bfS(\Qcoh\bbX)\ar[d]^-{\bfS f_*}\ar[rr]^-I\ar[d]&&\bfK(\Inj\bbX)
\ar[rr]^-Q\ar[d]^{\hat\bfR f_*}&&\bfD(\Qcoh\bbX)\ar[d]^{\bfR f_*}\\
\bfS(\Qcoh\bbY)\ar[rr]^-I&&\bfK(\Inj\bbY)
\ar[rr]^-Q&&\bfD(\Qcoh\bbY)\\
\bfS(\Qcoh\bbY)\ar[d]\ar[d]^-{\bfS f^!}&&\bfK(\Inj\bbY)\ar[ll]_-{I_\r}
\ar[d]^{\widehat{f^!}}&&\bfD(\Qcoh\bbY)\ar[ll]_-{Q_\r}\ar[d]^{f^!}\\
\bfS(\Qcoh\bbX)&&\bfK(\Inj\bbX)\ar[ll]_-{I_\r}
&&\bfD(\Qcoh\bbX)\ar[ll]_-{Q_\r}}$$
\end{thm}

It seems an interesting project to study the functor $\bfS f_*$, for
instance to find out when it is an equivalence. The following result
demonstrates the geometric content of this question; it generalizes a
result of Orlov for the bounded stable derived category \cite{O}.

\begin{thm}
Let $\bbY$ be a seperated noetherian scheme of finite Krull
dimension. If $f\colon \bbX\to\bbY$ denotes the inclusion of an open
subscheme which contains all singular points of $\bbY$, then $\bfS
f_*\colon \bfS(\Qcoh\bbX)\to\bfS(\Qcoh\bbY)$ is an equivalence.
\end{thm}

Despite the title of this paper and the algebraic geometric
formulation of the main results, there is another source of serious
interest in stable categories. Take a finite group $G$ and a field
$k$. A classical object in modular representation theory is the stable
module category $\uMod kG$ of the group algebra $kG$. We shall see
that this stable category is equivalent to the stable derived category
of the full module category $\Mod kG$. Using a slightly different
setting, Hovey, Palmieri, and Strickland studied the functor
$$I_\la\colon\bfK(\Inj kG)\lto \bfS(\Mod kG)\cong\uMod kG$$ in their
work on axiomatic stable homotopy theory \cite{HPS}.  Note that
$\bfK(\Inj kG)$ carries a commutative tensor product and the (graded)
endomorphism ring of its unit is nothing but the group cohomology ring
$H^*(G,k)$. Therefore, $\bfK(\Inj kG)$ seems to be the right object
for studying representations of $G$ via methods from commutative
algebra. In fact, the composite 
$$\bfD(\Mod kG)\xto{I_\la\comp Q_\r}\bfS(\Mod
kG)\xto{Z^0}\uMod kG$$ plays a crucial role in recent work of Benson
and Greenlees \cite{BG}.

Having stated some of the main results, let us sketch the outline of
this paper.  The paper deals with locally noetherian Grothendieck
categories and covers therefore various applications, for instance in
algebraic geometry or representation theory. Thus we fix a locally
noetherian Grothendieck category $\A$ and study the recollement
\begin{equation}\label{eq:recol}
\xymatrix{\bfS(\A)\ar[r]&\bfK(\Inj\A)
\ar[r]\ar@<.75ex>[l]\ar@<-.75ex>[l]&\bfD(\A).
\ar@<.75ex>[l]\ar@<-.75ex>[l]}
\end{equation}
More specifically, we begin in Section~\ref{se:htp} with the basic
properties of the homotopy category $\bfK(\Inj\A)$. The recollement
(\ref{eq:recol}) is established in Sections~\ref{se:loc} and
\ref{se:rec}. In Section~\ref{se:sta}, we discuss the essential
properties of the stable derived category $\bfS(\A)$. Then we extend
derived functors in Section~\ref{se:der}, and Section~\ref{se:gor} is
devoted to studying Gorenstein injective approximations and Tate
cohomology. In the final Section~\ref{se:kG}, we indicate the
relevance of the stable derived category in modular representation
theory. There is an appendix which provides additional material about DG
categories. Another appendix discusses homotopically minimal
complexes.

\section{The homotopy category of injectives}\label{se:htp}

We fix a {\em locally noetherian} Grothendieck category $\A$. Thus
$\A$ is an abelian Grothendieck category and has a set $\A_0$ of
noetherian objects which {\em generate} $\A$, that is, every object in
$\A$ is a quotient of a coproduct of objects in $\A_0$.  We denote by
$\noeth\A$ the full subcategory formed by the noetherian objects in
$\A$, and $\Inj\A$ denotes the full subcategory of injective
objects. Note that $\Inj\A$ is closed under taking coproducts.

We write $\bfK(\A)$ for the homotopy category and $\bfD(\A)$ for the
derived category of unbounded complexes in $\A$; for their definitions
and basic properties, we refer to \cite{V}.  We do
not distinguish between an object in $\A$ and the corresponding
complex concentrated in degree zero in the homotopy category
$\bfK(\A)$. The inclusion $\noeth\A\to\A$ induces  a fully faithful functor
$$\bfD^b(\noeth\A)\lto\bfD(\A)$$ which identifies $\bfD^b(\noeth\A)$
with the full subcategory of objects $X$ in $\bfD(\A)$ such that
$H^nX$ is noetherian for all $n$ and $H^nX=0$ for almost all
$n\in\bbZ$; see \cite[Proposition~III.2.4.1]{V}.

In this section, we study the basic properties of the homotopy category
$\bfK(\Inj\A)$. We shall see that this category solves a completion
problem for the triangulated category $\bfD^b(\noeth\A)$.  Let us
begin with some elementary observations.

\begin{lem}\label{le:com} 
Let $A$ be an object in $\A$ and denote by $iA$ an injective
resolution. Then the natural map
\begin{equation}\label{eq:iX}
\Hom_{\bfK(\A)}(iA,X)\lto\Hom_{\bfK(\A)}(A,X)
\end{equation} 
is an isomorphism for all $X$ in $\bfK(\Inj\A)$.  Therefore $iA$ is a
compact object in $\bfK(\Inj\A)$ if $A$ is noetherian.
\end{lem}
\begin{proof} 
Denote for any $n\in\mathbb Z$ by $\s^{\geqslant n}X$ the truncation satisfying
\begin{eqnarray*}
(\s^{\geqslant n}X)^p=
\begin{cases}X^p &\textrm{if }p\geqslant n,\\
0&\textrm{if } p<n.\\
\end{cases}
\end{eqnarray*}
We complete the map $A\to iA$ to an exact triangle
$$aA\lto A\lto iA\lto\Si(aA)$$ and obtain
$$\Hom_{\bfK(\A)}(aA,X)\cong\Hom_{\bfK(\A)}(aA,\s^{\geqslant -1}X)=0$$
since $aA$ is acyclic and concentrated in non-negative degrees. Thus
$\Hom_{\bfK(\A)}(iA,X)\cong\Hom_{\bfK(\A)}(A,X)$.

Now assume that $A$ is noetherian. Clearly, $A$ is a compact
object in $\A$ and therefore a compact object in $\bfK(\A)$. The
isomorphism (\ref{eq:iX}) shows that $iA$ is a compact object in
$\bfK(\Inj\A)$.
\end{proof}

\begin{lem}\label{le:gen}
Let $X$ be a non-zero object in $\bfK(\Inj\A)$. 
Then there exists a noetherian object $A$ in $\A$ such that
$\Hom_{\bfK(\A)}(A,\Si^nX)\neq 0$ for some $n\in\mathbb Z$.
\end{lem}
\begin{proof}
Suppose first $H^nX\neq 0$ for some $n$.  Choose a noetherian object
$A$ and a map $A\to Z^nX$ inducing a non-zero map $A\to H^nX$. We
obtain a chain map $A\to \Si^nX$ which induces a non-zero element
in $\Hom_{\bfK(\A)}(A,\Si^nX)$.

Now suppose $H^nX=0$ for all $n$.  We can choose $n$ such that $Z^nX$
is non-injective.  Using Baer's criterion, there exists a noetherian
object $A$ in $\A$ such that $\Ext^1_\A(A,Z^nX)$ is non-zero.  Now
observe that
$$\Hom_{\bfK(\A)}(A,\Si^{n+p}X)\cong\Ext^p_\A(A,Z^nX)$$ for all
$p\geqslant 1$.  Thus $\Hom_{\bfK(\A)}(A,\Si^{n+1}X)\neq 0$.  This completes
the proof.
\end{proof}

Let $\T$ be a triangulated category with arbitrary coproducts.  Recall
that an object $X$ in $\T$ is {\em compact} if $\Hom_\T(X,-)$
preserves all coproducts. The triangulated category is {\em compactly
generated} if there is a set $\T_0$ of compact objects such that
$\Hom_\T(X,\Si^nY)=0$ for all $X\in\T_0$ and $n\in\mathbb Z$ implies
$Y=0$ for every object $Y$ in $\T$.

\begin{prop}\label{pr:gen}
Let $\A$ be a locally noetherian Grothendieck category, and let $\bfK^c(\Inj\A)$
denote the full subcategory of compact objects in $\bfK(\Inj\A)$. 
\begin{enumerate}
\item The triangulated category $\bfK(\Inj\A)$ is compactly generated.
\item The canonical functor $\bfK(\A)\to\bfD(\A)$ induces an
equivalence $$\bfK^c(\Inj\A)\stackrel{\sim}\lto\bfD^b(\noeth\A).$$
\end{enumerate}
\end{prop}
\begin{proof} 
It follows from Lemma~\ref{le:com} and Lemma~\ref{le:gen} that
$\bfK(\Inj\A)$ is compactly generated. A standard argument shows that
$\bfK^c(\Inj\A)$ equals the thick subcategory of $\bfK(\Inj\A)$ which
is generated by the injective resolutions of the noetherian objects in
$\A$; see \cite[2.2]{N}.  The equivalence
$\bfK^+(\Inj\A)\to\bfD^+(\A)$ restricts to an equivalence
$\bfK^{+,b}(\Inj\A)\to\bfD^b(\A)$ and identifies $\bfK^c(\Inj\A)$ with
$\bfD^b(\noeth\A)$.
\end{proof}

Note that we obtain a functor $\bfD^b(\noeth\A)\to\bfK(\Inj\A)$ which
identifies $\bfD^b(\noeth\A)$ with the full subcategory of compact
objects.  Therefore the formation of the category $\bfK(\Inj\A)$
solves a completion problem which we explain by an analogy. The
category $\A$ is a completion of $\noeth\A$ in the following sense.
\begin{itemize}
\item $\A$ is an additive category with filtered colimits.
\item The inclusion $\noeth\A\to\A$ identifies $\noeth\A$ with the
  full subcategory of finitely presented objects.
\item $\A$ coincides with the smallest subcategory which contains all
  finitely presented objects and is closed under forming filtered
  colimits.
\end{itemize}
Recall that an object $X$ in $\A$ is {\em finitely presented} if the
functor $\Hom_\A(X,-)$ preserves filtered colimits. Similarly, 
we have the following for $\T=\bfK(\Inj\A)$.
\begin{itemize}
\item $\T$ is a triangulated category with coproducts.
\item The functor $\bfD^b(\noeth\A)\to\T$ identifies
  $\bfD^b(\noeth\A)$ with the full subcategory of compact objects.
\item $\T$ coincides with the smallest subcategory which contains all
  compact objects and is closed under forming triangles and coproducts.
\end{itemize}
The category $\A$ is, up to an equivalence, uniquely determined by
$\noeth\A$. It would be interesting to know to what extent
$\bfK(\Inj\A)$ is uniquely determined by $\bfD^b(\noeth\A)$.

\begin{exm}
Suppose there is a noetherian object $A$ in $\A$ such that
$\bfD^b(\noeth\A)$ is generated by $A$, that is, there is no proper
thick subcategory containing $A$. Take an injective resolution $iA$
and denote by $\END_{\A}(A)$ the endomorphism DG algebra of $iA$. Then
$\HOM_\A(iA,-)$ induces an equivalence between $\bfK(\Inj\A)$ and the
derived category $\bfD_\dg(\END_\A(A))$ of DG $\END_{\A}(A)$-modules;
see \cite{Ke}. If one replaces a single generator by a set of
generating objects, then one obtains an analogue which involves a DG
category instead of a DG algebra. In particular, $\noeth\A$ carries
the structure of a DG category such that $\bfK(\Inj\A)$ and
$\bfD_\dg(\noeth\A)$ are equivalent. We refer to Appendix~\ref{ap:dg} for
details.
\end{exm}

\begin{exm} 
Let $G$ be a finite $p$-group and $k$ be a field of characteristc
$p>0$.  We consider the category $\A=\Mod kG$ of modules over the
group algebra $kG$.  Take an injective resolution $ik$ of the trivial
representation $k$, and denote by $\END_{kG}(k)$ the endomorphism DG
algebra of $ik$. Then its derived category $\bfD_\dg(\END_{kG}(k))$ is
equivalent to $\bfK(\Inj\A)$. The tensor product $\otimes_k$ on $\A$
restricts to a product on $\Inj\A$ and induces therefore a (total)
tensor product on $\bfK(\Inj\A)$.  On the other hand, the
$E_\infty$-structure of $\END_{kG}(k)$ induces a product on
$\bfD_\dg(\END_{kG}(k))$. We conjecture that these products are naturally
isomorphic.
\end{exm}

\begin{exm}
Let $\La$ be a finite dimensional algebra over a field $k$. Then
$E=\Hom_k(\La^\op,k)$ is an injective cogenerator for $\A=\Mod\La$, and
$\Hom_\La(E,-)$ induces an equivalence $\Inj\A\to\Proj\A$ since
$\Hom_\La(E,E)\cong\La$. Thus the homotopy category $\bfK(\Proj\A)$ is
compactly generated. For more on  $\bfK(\Proj\A)$, see \cite{J,J2}.
\end{exm}

\section{A localization sequence}\label{se:loc}

Let $\A$ be a locally noetherian Grothendieck category and let
$$\bfK_\ac(\Inj\A)=\bfK(\Inj\A)\cap\bfK_\ac(\A),$$ where
$\bfK_\ac(\A)$ denotes the full subcategory formed by all acyclic
complexes in $\bfK(\A)$.  In this section, we prove that the canonical
functors
$$I\colon\bfK_\ac(\Inj\A)\xto{\inc}\bfK(\Inj\A) \quad\textrm{and}\quad
Q\colon\bfK(\Inj\A)\xto{\inc} \bfK(\A)\xto{\can}\bfD(\A)$$ form a
localization sequence
\begin{equation}\label{eq:locseq}
\xymatrix{\bfK_\ac(\Inj\A)\ar[r]^-I&\bfK(\Inj\A)\ar[r]^-Q&\bfD(\A).}
\end{equation}
Let us start with some preparations. In particular, we need to give
the definition of a localization sequence.

\begin{defn}\label{de:locseq}
We say that a sequence
$$\xymatrix{\T'\ar[r]^-F&\T\ar[r]^-G&\T''}$$ of exact functors
between triangulated categories is a {\em localization sequence} if
the following holds.
\begin{enumerate}
\item[(L1)] The functor $F$ has a right adjoint $F_\r\colon\T\to\T'$
satisfying $F_\r\comp F\cong\Id_{\T'}$.
\item[(L2)] The functor $G$ has a right adjoint $G_\r\colon\T''\to\T$
satisfying $G\comp G_\r\cong\Id_{\T''}$.
\item[(L3)] Let $X$ be an object in $\T$. Then $GX=0$ if and only if
$X\cong FX'$ for some $X'\in\T'$.
\end{enumerate}
The sequence $(F,G)$ of functors  is called {\em colocalization
sequence} if the sequence $(F^\op,G^\op)$ of opposite functors is a
localization sequence.
\end{defn}

The basic properties of a localization sequence are the following
\cite[II.2]{V}.
\begin{enumerate}
\item The functors $F$ and $G_\r$ are fully faithful.
\item Identify $\T'=\Im F$ and $\T''=\Im G_\r$. Given objects $X,Y\in\T$, then 
\begin{eqnarray*}
X\in\T'&\iff&\Hom_\T(X,\T'')=0,\\ 
Y\in\T''&\iff&\Hom_\T(\T',Y)=0.
\end{eqnarray*}
\item Identify $\T'=\Im F$. Then the functor $G$ induces an equivalence
$\T/\T'\to\T''$.
\item  Let $X$ be an object in $\T$. Then there is an exact triangle
$$(F\comp F_\r)X\lto X\lto (G_\r\comp G)X\lto\Si((F\comp F_\r)X)$$
which is functorial in $X$.
\item The sequence
$$\xymatrix{\T''\ar[r]^-{G_\r}&\T\ar[r]^-{F_\r}&\T'}$$
is a colocalization sequence.
\end{enumerate}

The next lemma is well known; it provides useful criteria for a
sequence to be a localization sequence. Recall that a full subcategory
of a triangulated category is {\em thick} if it is a triangulated
subcategory which is closed under taking direct factors.

\begin{lem}\label{le:locseq}
Let $\T$ be a triangulated category and $\S$ be a thick
subcategory. Then the following are equivalent.
\begin{enumerate}
\item The sequence $\S\xto{\inc}\T\xto{\can}\T/\S$ is a localization
sequence.
\item The inclusion functor $\S\to\T$ has a right adjoint.
\item The quotient functor $\T\to\T/\S$ has a right adjoint.
\end{enumerate}
\end{lem}
\begin{proof}
Condition (1) implies (2) and (3). Also, (2) and (3) together imply (1). Thus
we need to show that (2) and (3) are equivalent. Let us write
$F\colon\S\to\T$ and $G\colon\T\to\T/\S$ for the functors which are
involved.

(2) $\Rightarrow$ (3): We obtain a functor $L\colon\T\to\T$ by
completing for each $X$ in $\T$ the natural map $(F\comp F_\r) X\to X$
to an exact triangle
$$(F\comp F_\r) X\lto X\lto LX\lto\Si((F\comp F_\r) X).$$ The functor $L$
annihilates $\S$ and factors therefore through $G$ via an exact
functor $G_\r\colon\T/\S\to\T$. This is a right adjoint of $G$. In
fact, for each pair of objects $X$ in $\T$ and $Y$ in $\T/\S$, the
natural map
$$\Hom_{\T/\S}(GX,Y)\lto\Hom_\T(LX,G_\r Y)\lto\Hom_\T(X,G_\r Y)$$
is bijective.

(3) $\Rightarrow$ (2): We obtain a right
adjoint $F_\r\colon\T\to\S$ for the inclusion $F$ by completing for each $X$ in
$\T$ the natural map $X\to (G_\r\comp G)X$ to an exact triangle
$$F_\r X\lto X\lto (G_\r\comp G)X\lto\Si(F_\r X).$$ Note that $F_\r X$
belongs to $\S$ since $G(F_\r X)=0$.
\end{proof}

We need to construct left and right adjoints for functors starting in
a compactly generated triangulated category. Our basic tool for this
is the following result which is due to Neeman.

\begin{prop}\label{pr:adjoints}
Let $F\colon\S\to\T$ be an exact functor between triangulated
categories, and suppose $\S$ is compactly generated.
\begin{enumerate}
\item There is a right adjoint $\T\to\S$ if and only if $F$ preserves
all coproducts.
\item There is a left adjoint $\T\to\S$ if and only if $F$ preserves
all products.
\end{enumerate}
\end{prop}
\begin{proof}
For (1), see \cite[Theorem~4.1]{N1}. The proof of (2) is analogous and
uses covariant Brown representability \cite[Theorem~8.6.1]{N2}; see
also \cite{K}.
\end{proof}

We record a similar result for later use.

\begin{prop}\label{pr:adjoints2}
Let $\T$ be a compactly generated triangulated category and $\S_0$ be
a set of objects in $\T$. Denote by $\U$ the full subcategory of
objects $Y$ in $\T$ such that $\Hom_\T(\Si^nX,Y)=0$ for all $X\in\S_0$
and $n\in\bbZ$. Then the inclusion $\U\to\T$ has a left adjoint.
\end{prop}
\begin{proof}
The localizing subcategory $\S$ generated by $\S_0$ is well generated
and the inclusion $\S\to\T$ has therefore a right adjoint; see
\cite{N2}. We obtain a localization sequence
$\S\xto{\inc}\T\xto{\can}\T/\S$ by Lemma~\ref{le:locseq}, and the
right adjoint of the canonical functor $\T\to\T/\S$ identifies $\T/\S$
with $\U$.
\end{proof}

There is a useful criterion when a left adjoint preserves compactness.

\begin{lem}\label{le:adjoint}
Let $F\colon\S\to\T$ be an exact functor between compactly generated
triangulated categories which has a right adjoint $G$.  Then $F$
preserves compactness if and only if $G$ preserves coproducts.
\end{lem}
\begin{proof}
  See \cite[Theorem~5.1]{N1}. 
\end{proof}

The following result establishes the localization sequence for the
homotopy category of injective objects.

\begin{prop}\label{pr:seq}
Let $\A$ be a locally noetherian Grothendieck category.  
Then the canonical functors $\bfK_\ac(\Inj\A)\to\bfK(\Inj\A)$ and
$\bfK(\Inj\A)\to\bfD(\A)$ form a localization sequence
$$
\xymatrix{\bfK_\ac(\Inj\A)\ar[r]^-I&\bfK(\Inj\A)\ar[r]^-Q&\bfD(\A).}
$$
\end{prop}
\begin{proof}
We know from Proposition~\ref{pr:gen} that $\bfK(\Inj\A)$ is compactly
generated.  In addition, we use Lemma~\ref{le:locseq} and
Proposition~\ref{pr:adjoints}.  The inclusion
$J\colon\bfK(\Inj\A)\to\bfK(\A)$ preserves products and has therefore
a left adjoint $J_\la$ satisfying $J_\la\comp
J\cong\Id_{\bfK(\Inj\A)}$.  We obtain a localization sequence
$$
\xymatrix{\K\ar[r]^-\inc&\bfK(\A)\ar[r]^-{J_\la}&\bfK(\Inj\A)}
$$ where $\K$ denotes the kernel of $J_\la$. Thus
$$\Hom_{\bfK(\A)}(X,Y)=0\quad\textrm{for all}\quad
X\in\K\quad\textrm{and}\quad Y\in\bfK(\Inj\A).$$ This implies
$\K\subseteq \bfK_\ac(\A)$ and gives the following commutative
diagram of exact functors.
$$\xymatrix{\K\ar[d]^\inc\ar[rr]^-{\inc}&&
\bfK(\A)\ar@{=}[d]\ar[rr]^-{J_\la}&&\bfK(\Inj\A)\ar[d]^F\\
\bfK_\ac(\A)\ar[rr]^-{\inc}&&\bfK(\A)\ar[rr]^-\can&&\bfD(\A)}
$$ The functor $F$ is induced by the canonical functor
$\bfK(\A)\to\bfD(\A)$, and we have $F\cong Q$ since $J_\la\comp
J\cong\Id_{\bfK(\Inj\A)}$. Moreover, $F$ preserves coproducts and has
therefore a right adjoint $F_\r$. The composite $J\comp F_\r$ is a
right adjoint for the canonical functor $\bfK(\A)\to\bfD(\A)$. This
implies $F\comp F_\r\cong\Id_{\bfD(\A)}$. On the other hand,
$\bfK_\ac(\Inj\A)$ is the kernel of $F$. Thus we conclude that the
sequence (\ref{eq:locseq}) is a localization sequence.
\end{proof}

We add some useful remarks which are immediate consequences.

\begin{rem}\label{re:quotient}
  Let $J_\la\colon\bfK(\A)\to\bfK(\Inj\A)$ be the left adjoint of the
  inclusion $\bfK(\Inj\A)\to\bfK(\A)$. Then the composite $Q\comp
  J_\la$ is naturally isomorphic to the canonical functor
  $\bfK(\A)\to\bfD(\A)$.
\end{rem}

\begin{rem}\label{re:inverse}
The right adjoint $Q_\r$ of $Q$ induces an equivalence
$$\bfD^b(\noeth\A)\stackrel{\sim}\lto\bfK^c(\Inj\A)$$ which is a
quasi-inverse for the equivalence
${\bfK^c(\Inj\A)}\to\bfD^b(\noeth\A)$ induced by $Q$.
\end{rem}

Let us denote by $\bfK_\inj(\A)$ the full subcategory of complexes $Y$
in $\bfK(\Inj\A)$ such that $\Hom_{\bfK(\A)}(X,Y)=0$ for all acyclic
complexes $X$ in $\bfK(\A)$.  Following Spaltenstein's terminology
\cite{S}, the objects in $\bfK_\inj(\A)$ are precisely the {\em
K-injective} complexes having injective components. There are various
results about K-injective resolutions in the literature; see for
instance \cite{S,BN}. The following is certainly not the most general
one; however it is sufficient in our context.

\begin{cor}\label{co:i}
The inclusion $\bfK_\inj(\A)\to\bfK(\A)$ has a left adjoint
$i\colon\bfK(\A)\to\bfK_\inj(\A)$ which has the following properties.
\begin{enumerate}
\item Every object $X$ in $\bfK(\A)$ fits into an exact triangle 
$$aX\lto X\lto iX\lto\Si(aX)$$ such that $aX$ is an acyclic complex.
\item The functor $i\colon\bfK(\A)\to\bfK_\inj(\A)$ 
induces an equivalence
$$\bfD(\A)=\bfK(\A)/\bfK_\ac(\A)\stackrel{\sim}\lto\bfK_\inj(\A).$$
\item We have for all $X,Y$ in $\bfK(\A)$
$$\Hom_{\bfD(\A)}(X,Y)\cong\Hom_{\bfK(\A)}(X,iY).$$ 
\end{enumerate}
\end{cor}
\begin{proof}
Put $iX=Q_\r X$ for each $X$ in $\bfK(\A)$, where $Q_\r$ denotes the
right adjoint of $Q\colon \bfK(\Inj\A)\to\bfD(\A)$. The properties of
the functor $i$ follow from the fact that $J\comp Q_\r$ is a right
adjoint of the canonical functor $\bfK(\A)\to\bfD(\A)$. In particular,
we see that $iX$ is a K-injective complex.
\end{proof}

The functor
$$R\colon\bfD(\A)=\bfK(\A)/\bfK_\ac(\A)\stackrel{\sim}
\lto\bfK_\inj(\A)\stackrel{\inc}\lto\bfK(\A)$$ provides a right
adjoint for the canonical functor  $\bfK(\A)\to\bfD(\A)$.
Let us mention as an application that the right derived functor
of any additive functor $F\colon\A\to\B$ is obtained as composite
$$\bfR F\colon \bfD(\A)\stackrel{R}\lto\bfK(\A)\stackrel{\bfK(F)}
\lto\bfK(\B)\stackrel{\can}\lto\bfD(\B).$$

\begin{exm}
Suppose every object in $\A$ has finite injective dimension. Then the
functor $\bfK(\Inj\A)\to\bfD(\A)$ is an equivalence since
$\bfK_\ac(\Inj\A)=0$.  In particular, the compact objects in
$\bfD(\A)$ are precisely those from $\bfD^b(\noeth\A)$.
\end{exm}

\begin{exm}
Suppose products in $\A$ are exact. For instance, let $\A$ be a module
category.  Then one can show that $\bfK_\inj(\A)$ is the smallest
triangulated subcategory of $\bfK(\A)$ which is closed under taking
products and contains the injective objects of $\A$ (viewed as
complexes concentrated in degree zero).
\end{exm}

\section{A recollement}\label{se:rec}
In this section, we provide a criterion for $\A$ such that the sequence
\begin{equation*}
\xymatrix{\bfK_\ac(\Inj\A)\ar[r]^-I&\bfK(\Inj\A)\ar[r]^-Q&\bfD(\A)}
\end{equation*}
induces  a recollement 
$$\xymatrix{\bfK_\ac(\Inj\A)\ar[r]&\bfK(\Inj\A)
\ar[r]\ar@<.75ex>[l]\ar@<-.75ex>[l]&\bfD(\A)
\ar@<.75ex>[l]\ar@<-.75ex>[l]}$$ in the sense of \cite{BBD}. It is
important to note that one cannot expect a recollement
\begin{equation}\label{eq:reco}
\xymatrix{\bfK_\ac(\A)\ar[r]&\bfK(\A)
\ar[r]\ar@<.75ex>[l]\ar@<-.75ex>[l]&\bfD(\A)
\ar@<.75ex>[l]\ar@<-.75ex>[l]}
\end{equation}
without severe restrictions on $\A$; see Example~\ref{ex:exact}. In fact, a recollement
(\ref{eq:reco}) implies that a product of exact sequences in $\A$ remains exact.

We begin with a lemma.

\begin{lem}\label{le:der_comp}
Let $\A$ be a locally noetherian Grothendieck category. Then a compact
object in $\bfD(\A)$ belongs to $\bfD^b(\noeth\A)$.
\end{lem}
\begin{proof}
Suppose $X$ is compact in $\bfD(\A)$. We need to show that $H^nX$ is
noetherian for all $n$, and that $H^nX$ vanishes for almost all $n$ in
$\bbZ$. We have for any injective object $E$ in $\A$ an isomorphism
$$\Hom_{\bfD(\A)}(X,E)\cong\Hom_\A(H^0X,E).$$ Therefore
$\Hom_\A(H^0X,-)$ preserves coproducts in $\Inj\A$. This implies that
each $H^nX$ is noetherian; see \cite{Re}. Now fix for each $n$ an
injective envelope $H^nX\to E(H^nX)$ and consider the induced map
$$\a\colon X\lto\prod_{n\in\bbZ}\Si^{-n}E(H^nX)$$ in $\bfD(\A)$. The
canonical map
$$\coprod_{n\in\bbZ}\Si^{-n}E(H^nX)\lto\prod_{n\in\bbZ}\Si^{-n}E(H^nX)$$
is an isomorphism in $\bfD(\A)$, and therefore $\a$ factors though a
finite number of factors in $$\prod_{n\in\bbZ}\Si^{-n}E(H^nX).$$ Thus
$H^nX$ vanishes for almost all $n$ in $\bbZ$, and the proof is
complete.
\end{proof}

We denote by $\bfD^c(\A)$ the full subcategory of $\bfD(\A)$ which is
formed by all compact objects.

\begin{thm}\label{th:products}
Let $\A$ be a locally noetherian Grothendieck category and suppose
$\bfD(\A)$ is compactly generated. Then the canonical functor $Q\colon
\bfK(\Inj\A)\to\bfD(\A)$ has a left adjoint and therefore the sequence
\begin{equation*}
\xymatrix{\bfK_\ac(\Inj\A)\ar[r]^-I&\bfK(\Inj\A)\ar[r]^-Q&\bfD(\A)}
\end{equation*}
is a colocalization sequence.
\end{thm}
\begin{proof}
Let $\K$ be the localizing subcategory of $\bfK(\Inj\A)$ which is
generated by all compact objects $X$ in $\bfK(\Inj\A)$ such that $QX$
is compact in $\bfD(\A)$. We claim that $Q|_\K\colon\K\to\bfD(\A)$ is
an equivalence. First note that $\K$ and $\bfD(\A)$ are both compactly
generated. We have seen in Lemma~\ref{le:der_comp} that
$$\bfD^c(\A)\subseteq\bfD^b(\noeth\A),$$ 
and $Q$ induces an equivalence
$$\bfK^c(\Inj\A)\stackrel{\sim}\lto\bfD^b(\noeth\A),$$ by
Proposition~\ref{pr:gen}. Thus $Q$ induces an equivalence between the
subcategories of compact objects in $\K$ and $\bfD(\A)$. Then a
standard argument shows that $Q|_\K$ is an equivalence since $Q$
preserves all coproducts. Now fix a left adjoint
$L\colon\bfD(\A)\to\K$. We claim that the composite
$$\bfD(\A)\stackrel{L}\lto \K\stackrel{\inc}\lto \bfK(\Inj\A)$$ is a
left adjoint for $Q$. To see this, consider for objects $X$ in $\bfD(\A)$
and $Y$ in $\bfK(\Inj\A)$ the natural map
$$\a_{X,Y}\colon\Hom_{\bfK(\Inj\A)}(LX,Y)\lto
\Hom_{\bfD(\A)}(QLX,QY)\stackrel{\sim}\lto\Hom_{\bfD(\A)}(X,QY)$$
which is induced by $Q$. If $X$ and $Y$ are compact, then $\a_{X,Y}$
is bijective by Proposition~\ref{pr:gen}. We use a standard argument
to show that $\a_{X,Y}$ is bijective for arbitrary $X$ and $Y$. Fix a
compact object $X$. Then the objects $Y$ such that $\a_{X,Y}$ is
bijective form a triangulated subcategory which is closed under taking
coproducts and contains all compact objects. Thus $\a_{X,Y}$ is
bijective for all $Y$ because $\bfK(\Inj\A)$ is compactly generated.
Now fix any object $Y$. The same argument shows that $\a_{X,Y}$ is
bijective for all $X$ because $\bfD(\A)$ is compactly generated. We
conclude that $Q$ has a left adjoint. Moreover, Lemma~\ref{le:locseq}
implies that $I$ and $Q$ form a colocalization sequence.
\end{proof}


Following Beilinson, Bernstein, and Deligne \cite{BBD}, we say that a
sequence
\begin{equation}\label{eq:rec}
\xymatrix{\T'\ar[r]&\T\ar[r]&\T''}
\end{equation}
of exact functors between triangulated categories induces a {\em
recollement}
$$\xymatrix{\T'\ar[r]&\T \ar[r]\ar@<.75ex>[l]\ar@<-.75ex>[l]&\T''
\ar@<.75ex>[l]\ar@<-.75ex>[l]}$$ if the sequence (\ref{eq:rec}) is a
localization sequence and a colocalization sequence in the sense of
Definition~\ref{de:locseq}.

\begin{cor} 
Let $\A$ be a locally noetherian Grothendieck category and suppose
$\bfD(\A)$ is compactly generated. Then the sequence
$$\xymatrix{\bfK_\ac(\Inj\A)\ar[r]^-I&\bfK(\Inj\A)\ar[r]^-Q&\bfD(\A)}$$
induces a recollement
$$\xymatrix{\bfK_\ac(\Inj\A)\ar[r]&\bfK(\Inj\A)
\ar[r]\ar@<.75ex>[l]\ar@<-.75ex>[l]&\bfD(\A).
\ar@<.75ex>[l]\ar@<-.75ex>[l]}$$
\end{cor}

\begin{cor}\label{co:products}
Let $\A$ be a locally noetherian Grothendieck category and suppose
$\bfD(\A)$ is compactly generated. Then a product of acyclic complexes
of injective objects in $\A$ is acyclic.
\end{cor}

Let us give a criterion for $\A$ such that the derived category
$\bfD(\A)$ is compactly generated.

\begin{lem}
Let $\A$ be a locally noetherian Grothendieck category. Suppose there
is a set $\A_0$ of objects in $\A$ which are compact when viewed as
objects in $\bfD(\A)$. If $\A_0$ generates $\A$, then $\bfD(\A)$ is
compactly generated by $\A_0$.
\end{lem}

The lemma is an immediate consequence of the following statement.

\begin{lem}\label{le:compgen}
Let $\A$ be a locally noetherian Grothendieck category and fix a set
$\A_0$ of generating objects. Let $X$ be a complex in $\A$ such that
$H^0X\neq 0$. Then there exists some object $A$ in $\A_0$ such that
$$\Hom_{\bfK(\A)}(A,X)\neq
0\quad\textrm{and}\quad\Hom_{\bfD(\A)}(A,X)\neq 0.$$
\end{lem}
\begin{proof}
Choose $A$ in $\A_0$ and a map $A\to Z^0X$ such that the composite
with $Z^0X\to H^0X$ is non-zero. This induces a non-zero element in
$$H^0(\Hom_\A(A,X))\cong\Hom_{\bfK(\A)}(A,X).$$ The second assertion
follows from the first since for any object $A$ in $\A$ we have
$$\Hom_{\bfD(\A)}(A,X)\cong\Hom_{\bfK(\A)}(A,iX)$$ and $H^0(iX)\cong H^0X$.
\end{proof}

We give examples of Grothendieck categories such that objects in
$\A$ become compact objects in $\bfD(\A)$.

\begin{exm}
  Let $\La$ be an associative ring. Denote by $\A=\Mod\La$ the
  category of (right) $\La$-modules and by $\proj\La$ the full
  subcategory of finitely generated projective $\La$-modules. Then
  every object in $\proj\La$ is compact when viewed as object in
  $\bfD(\A)$. Thus the inclusion $\bfD^b(\proj\La)\to\bfD(\A)$
  identifies $\bfD^b(\proj\La)$ with the full subcategory of compact
  objects in $\bfD(\A)$.  Suppose now that $\La$ is right
  noetherian.  Then the fully faithful functor $Q_\la\colon
  \bfD(\A)\to \bfK(\Inj\A)$ identifies $\bfD(\A)$ with the localizing
  subcategory of $\bfK(\Inj\A)$ which is generated by the injective
  resolution $i\La$ of $\La$.
\end{exm}

Let us return to the completion problem for triangulated categories
which has been addressed in Section~\ref{se:htp}. Keeping the analogy
between the completion with respect to filtered colimits and the
completion with respect to triangles und coproducts, we obtain the
following diagram for a right noetherian ring $\La$. The vertical
arrows denote completions and the horizontal ones the appropriate inclusions.
$$\xymatrix{ \proj\La\ar[d]\ar[r]&\mod\La\ar[d]&
  \bfD^b(\proj\La)\ar[d]\ar[r]&\bfD^b(\mod\La)\ar[d]\\
  \Flat\La\ar[r]&\Mod\La&\bfD(\Mod\La)\ar[r]^{Q_\la}&\bfK(\Inj\La) }$$
Here, $\Flat\La$ denotes the full subcategory of flat $\La$-modules,
which is the closure of $\proj\La$ under forming filtered colimits. 

\begin{exm}
Let $\bbX$ be a quasi-compact and separated scheme, and let $L$ be a
locally free sheaf of finite rank.  Then
$$\Hom_{\bfD(\Qcoh\bbX)}(L,-)\cong \Hom_{\bfD(\Qcoh\bbX)}(\mathcal
O_\bbX,L^\vee\otimes_{\mathcal O_\bbX}-) \cong
H^0(L^\vee\otimes_{\mathcal O_\bbX}-),$$ where $L^\vee=\HOM_{\mathcal
O_\bbX}(L,\mathcal O_\bbX)$. Thus $L$ is a compact object in
$\bfD(\Qcoh\bbX)$; see \cite{N1}. If $\bbX$ has an ample family of
line bundles, then the locally free sheaves of finite rank generate
$\Qcoh\bbX$.
\end{exm}

It would be interesting to know in which generality products of
acyclic complexes of injectives are acyclic. In fact, I do not know an
example where this property fails. However, it is important to
restrict to complexes of injectives. In order to illustrate this
point, let us include an example which shows that products in
$\Qcoh\bbX$ need not to be exact. I learned this example from Bernhard
Keller.

\begin{exm}\label{ex:exact}
  Let $k$ be a field and $\bbX=\bbP^1_k$ the projective line with
  homogeneous coordinate ring $S=k[x_0,x_1]$. For each $n\geqslant 0$, we
  have a canonical map
  $$\pi_n\colon\OO(-n)\otimes_k\Hom_\bbX(\OO(-n),\OO)\lto\OO$$
which is an epimorphism in $\Qcoh\bbX$.  We claim that the product
 $$\pi\colon\prod_{n\geqslant
 0}\bigl(\OO(-n)\otimes_k\Hom_\bbX(\OO(-n),\OO)\bigr)\lto\prod_{n\geqslant
 0}\OO$$ is not an epimorphism.  Taking graded global sections gives for
 each $n\geqslant 0$ the multiplication map
$$\Ga_*(\bbX,\pi_n)\colon S(-n)\otimes_k S_n\lto S$$ which is a map of
graded $S$-modules with cokernel of finite length. However, the
cokernel of
$$\Ga_*(\bbX,\pi)=\prod_{n\geqslant 0}\Ga_*(\bbX,\pi_n)$$ is not a torsion
module. The left adjoint of $\Ga_*(\bbX,-)$ is exact and takes
$\Ga_*(\bbX,\pi)$ to $\pi$. It follows that the cokernel of $\pi$ is
non-zero, because the left adjoint of $\Ga_*(\bbX,-)$ annihilates
exactly those $S$-modules which are torsion.
\end{exm}

\section{The stable derived category}\label{se:sta}

Let $\A$ be a locally noetherian Grothendieck category. We suppose
that $\bfD(\A)$ is compactly generated.

\begin{defn}
The {\em stable derived category} $\bfS(\A)$ of $\A$ is by definition
the full subcategory of $\bfK(\A)$ which is formed by all acyclic
complexes of injective objects in $\A$. The full subcategory of
compact objects is denoted by $\bfS^c(\A)$.
\end{defn}

In this section, we show that the stable derived category is compactly
generated, and the description of the category of compact objects
justifies our terminology.  Our basic tool is the (co)localization
sequence
$$\xymatrix{\bfS(\A)\ar[r]^-{I}& \bfK(\Inj\A)\ar[r]^-{Q}&\bfD(\A)}.$$
Thus we use the fact that $I$ and $Q$ have left adjoints $I_\la$,
$Q_\la$ and right adjoints $I_\r$, $Q_\r$.  The {\em stabilization
functor} is by definition the composite
$$S\colon\bfD(\A)\xto{I_\la\comp Q_\r}\bfS(\A).$$
We begin with the following lemma.

\begin{lem}\label{le:nat_transf}
  Let $\A$ be a locally noetherian Grothendieck category and suppose
  $\bfD(\A)$ is compactly generated. The functors
  $Q_\la,Q_\r\colon\bfD(\A)\to\bfK(\Inj\A)$ admit a natural
  transformation $\eta\colon Q_\la\to Q_\r$, and $\eta$ is an
  isomorphism when restricted to the subcategory of compact objects in
  $\bfD(\A)$.
\end{lem}
\begin{proof}
We have a natural isomorphism $\mu\colon
\Id_{\bfD(\A)}\stackrel{\sim}\to Q\comp Q_\la$.  The natural
transformation
$$Q_\la\comp Q\lto\Id_{\bfK(\Inj\A)}\lto Q_\r\comp Q$$ induces
for each $X$ in $\bfD(\A)$
a natural map 
$$\eta_X\colon Q_\la X\xto{Q_\la(\mu_X)} (Q_\la\comp Q)Q_\la X\lto
(Q_\r\comp Q)Q_\la X\xto{Q_\r(\mu^{-1}_X)} Q_\r X.$$ Note that $Q(\eta)$
induces an isomorphism
$$Q\comp Q_\la\stackrel{\sim}\lto Q\comp Q_\r.$$
We know from
Proposition~\ref{pr:gen} that $Q$ induces an equivalence
$$\bfK^c(\Inj\A)\stackrel{\sim}\lto\bfD^b(\noeth\A).$$ On the other
hand, $$Q_\la(\bfD^c(\A))\subseteq\bfK^c(\Inj\A)$$ since a left
adjoint preserves compactness if the right adjoint preserves
coproducts; see Lemma~\ref{le:adjoint}. Also,
$$Q_\r(\bfD^c(\A))\subseteq\bfK^c(\Inj\A),$$ since
$\bfD^c(\A)\subseteq\bfD^b(\noeth\A)$ by Lemma~\ref{le:der_comp}, and
$$Q_\r(\bfD^b(\noeth\A))=\bfK^c(\Inj\A)$$
by Remark~\ref{re:inverse}.
We conclude that $\eta|_{\bfD^c(\A)}$ is an isomorphism.
\end{proof}

\begin{prop}\label{pr:stable}
Let $\A$ be a locally noetherian Grothendieck category, and suppose
$\bfD(\A)$ is compactly generated.  Then we have a localization
sequence
\begin{equation}\label{eq:loc}
\xymatrix{\bfD(\A)\ar[r]^-{Q_\la}&
\bfK(\Inj\A)\ar[r]^-{I_\la}&\bfS(\A)}
\end{equation}
which induces the following commutative diagram.
$$\xymatrix{ \bfD^c(\A)\ar@{=}[d]\ar[rr]^-{\inc}&&\bfD^b(\noeth\A)
\ar[d]_-{\wr}^{Q_\r|_{\bfD^b(\noeth\A)}}
\ar[rr]^-{\can}&&\bfD^b(\noeth\A)/\bfD^c(\A)\ar[d]^-F\\
\bfD^c(\A)\ar[d]^-{\inc}\ar[rr]&&
\bfK^c(\Inj\A)\ar[d]^-{\inc}\ar[rr]&&\bfS^c(\A)\ar[d]^-{\inc}\\
{\bfD(\A)}\ar[rr]^-{Q_\la}&&\bfK(\Inj\A)\ar[rr]^-{I_\la}&&\bfS(\A)\\
}
$$
\end{prop}
\begin{proof}
It follows from Theorem~\ref{th:products} that the sequence
(\ref{eq:loc}) is a localization sequence.  Let us explain the
commutativity of the diagram. First observe that a left adjoint
preserves compactness if the right adjoint preserves coproducts; see
Lemma~\ref{le:adjoint}. Therefore $I_\la$ and $Q_\la$ preserve
compactness, and this explains the commutativity of the lower
squares. Now observe that
$$\bfD^c(\A)\subseteq\bfD^b(\noeth\A),$$
by Lemma~\ref{le:der_comp},
and that $Q_\r|_{\bfD^b(\noeth\A)}$ is a quasi-inverse for
$Q|_{\bfK^c(\Inj\A)}$.  It follows from Lemma~\ref{le:nat_transf} that
the upper left hand square commutes.  The functor $F$ is by definition
the unique functor making the upper right hand square commutative. It
exists because $I_\la\comp Q_\la=0$.
\end{proof}

We have seen that the stable derived category $\bfS(\A)$ is a
localization of the homotopy category $\bfK(\Inj\A)$. This has some
interesting consequences.

\begin{cor}\label{co:stable}
The stable derived category $\bfS(\A)$ is compactly generated, and the
functor $I_\la\comp Q_\r\colon\bfD(\A)\to\bfS(\A)$ induces (up to
direct factors) an equivalence
$$F\colon\bfD^b(\noeth\A)/\bfD^c(\A)\stackrel{\sim}\lto \bfS^c(\A).$$
\end{cor}
\begin{proof}
We know from Proposition~\ref{pr:gen} that $\bfK(\Inj\A)$ is compactly
generated. This property carries over to $\bfS(\A)$ since $I_\la$
sends a set of compact generators of $\bfK(\Inj\A)$ to a set of
compact generators of $\bfS(\A)$.  The functor $Q_\la$ identifies
$\bfD(\A)$ with the localizing subcategory of $\bfK(\Inj\A)$ which is
generated by all compact objects in the image of $Q_\la$. Now apply
the localization theorem of Neeman-Ravenel-Thomason-Trobaugh-Yao
\cite{N}. This result describes the category of compact objects of the
quotient $\bfS(\A)$, up to direct factors, as the quotient of the
compact objects in $\bfK(\Inj\A)$ modulo those from the localizing
subcategory. To be precise, $F$ is fully faithful and every object in
$\bfS^c(\A)$ is a direct factor of some object in the image of $F$.
\end{proof}

\begin{cor}\label{co:projs} The composite
$$\A\xto{\can}\bfD(\A)\xto{I_\la\comp Q_\r}\bfS(\A)$$
preserves all coproducts and annihilates the objects in $\A\cap\bfD^c(\A)$.
\end{cor}
\begin{proof}
The diagram in Proposition~\ref{pr:stable} shows that $I_\la\comp
Q_\r$ annihilates $\A\cap\bfD^c(\A)$. To show that $I_\la\comp Q_\r$
preserves all coproducts, observe that $Q_\r$ sends an object in $\A$
to an injective resolution. A coproduct of injective
resolutions is again an injective resolution, and the left adjoint $I_\la$
preserves all coproducts. This finishes the proof.
\end{proof}

Using the stabilization functor $S\colon\bfD(\A)\to\bfS(\A)$, we
define for objects $X,Y$ in $\bfD(\A)$ and $n\in\bbZ$ the {\em stable
cohomology group}
$$\uExt^n_\A(X,Y)=\Hom_{\bfK(\A)}(SX,\Si^n(SY)).$$
Note that in both
arguments each exact sequence in $\A$ induces a long exact sequence in
stable cohomology. We do not go into details but refer to our
discussion of Tate cohomology in Section~\ref{se:gor}.  In fact, both
cohomology theories coincide in case $\A$ satisfies some appropriate
Gorenstein property, and we shall see explicit formulae for the Tate
cohomology groups.

\begin{exm}\label{ex:pd}
  Suppose $\A$ is a module category. Then the stabilization functor
  annihilates all finitely generated projective modules, and all
  coproducts of such, by Corollary~\ref{co:projs}. Hence it annihilates
  all projective modules. Since $I_\la\comp Q_\r$ is an exact functor
  vanishing on projectives, it annihilates all bounded complexes of
  projective modules. In particular, all modules of finite projective
  dimension are annihilated. Similarly, if $\A$ is a category of
  quasi-coherent sheaves, then the stabilization functor annihilates
  all sheaves having a finite resolution with locally free sheaves.
\end{exm}

Given a noetherian scheme $\bbX$, the stable derived category
$\bfS(\Qcoh\bbX)$ vanishes if $\bbX$ is regular. Nonetheless, a
classical result of Bernstein-Gelfand-Gelfand \cite{BGG} shows that
stable derived categories are relevant when one studies regular
schemes. This is sketched in the following example.

\begin{exm}
  Let $\La$ be a Koszul algebra and $\La^!$ its Koszul dual. Then we
  have under appropriate assumptions an equivalence
  $\bfK(\Inj\La)\stackrel{\sim}\to\bfK(\Inj\La^!)$ which induces an
  equivalence $\bfD^b(\mod\La)\stackrel{\sim}\to\bfD^b(\mod\La^!)$
  when restricted to the full subcategories of compact objects
  \cite{BGS,Ke}. Note that we consider the categories of graded
  modules over $\La$ and $\La^!$ respectively. The classical example
  is the symmetric algebra $\La=SV$ of a $d$-dimensional space $V$
  over a field $k$, where $\La^!$ is the exterior algebra $\bigwedge
  V^*$ of the dual space $V^*$.  The equivalence
  $\bfK(\Inj\La)\stackrel{\sim}\to\bfK(\Inj\La^!)$ takes an injective
  resolution $ik$ of $\La_0=k$ to $\La^!$ and identifies the
  localizing subcategory $\K$ generated by $ik$ with the localizing
  subcategory generated by $\La^!$, which is $\bfD(\Mod\La^!)$. Note
  that the quotient $\bfK(\Inj\La)/\K$ identifies with the derived
  category of the quotient $\Mod\La/(\Mod\La)_0$, where $(\Mod\La)_0$
  denotes the subcategory of torsion modules. This quotient is
  equivalent to $\Qcoh\bbP_k^{d-1}$ by Serre's Theorem. Thus we obtain
  an equivalence
  $$\bfD(\Qcoh\bbP_k^{d-1})\stackrel{\sim}\lto
  \bfS(\Mod\mbox{$\bigwedge k^d$}).$$
  Note that $\bfS(\Mod\bigwedge
  k^d)$ is equivalent to the stable module category $\uMod\bigwedge
  k^d$ because the exterior algebra is self-injective; see
  Example~\ref{ex:selfinj}.  Passing to the subcategory of compact
  objects, one obtains the equivalence
$$\bfD^b(\coh\bbP_k^{d-1})\stackrel{\sim}\lto \umod \mbox{$\bigwedge
  k^d$}$$
of Bernstein-Gelfand-Gelfand \cite{BGG}, where
$\umod\bigwedge k^d$ denotes the stable category of all finite
dimensional $\bigwedge k^d$-modules.  This example generalizes to
non-commutative algebras, for instance, to Artin-Schelter
regular algebras \cite{J3}.
\end{exm}

\section{Extending derived functors}\label{se:der}
An additive functor $F\colon\A\to\B$ between locally noetherian
Grothendieck categories admits a right derived functor $\bfR
F\colon\bfD(\A)\to\bfD(\B)$. In this section, we extend this to a
functor $\hat\bfR F\colon \bfK(\Inj\A)\to\bfK(\Inj\B)$ and investigate
its right and left adjoints.  As an application, we consider for $F$
the direct image functor $f_*\colon \Qcoh\bbX\to\Qcoh\bbY$
corresponding to a morphism $f\colon\bbX\to\bbY$ between noetherian
schemes. We use the following functors
$$J\colon\bfK(\Inj\A)\stackrel{\inc}\lto\bfK(\A)\quad\textrm{and}\quad
Q\colon\bfK(\Inj\A)\stackrel{\inc}\lto\bfK(\A)\stackrel{\can}\lto\bfD(\A)$$
simultanously for $\A$ and $\B$. Moreover, we use the fact that both
functors have left and right adjoints. 

\begin{thm}\label{th:extensions}
Let $F\colon \A\to\B$ be an additive functor between locally
noetherian Grothendieck categories. Suppose $\bfD(\A)$ and $\bfD(\B)$
are compactly generated.  Then the composite
$$\hat\bfR
F\colon\bfK(\Inj\A)\stackrel{J}\lto\bfK(\A)\stackrel{\bfK(F)}\lto
\bfK(\B)\stackrel{J_\la}\lto\bfK(\Inj\B)$$ makes the following diagram
commutative.
$$\xymatrix{\bfD(\A)\ar[d]^{Q_\r}\ar[rr]^-{\bfR F}&&\bfD(\B)\\
\bfK(\Inj\A)\ar[rr]^-{\hat\bfR F}&&\bfK(\Inj\B)\ar[u]_Q}$$

\begin{enumerate} 
\item Suppose $F$ preserves coproducts. Then $\hat\bfR F$ preserves
coproducts and has therefore a right adjoint $(\hat\bfR F)_\r$.
\item Suppose $F$ and $\bfR F$ preserve coproducts. Then $\bfR F$ has
a right adjoint $(\bfR F)_\r$ making the following diagram commutative.
$$\xymatrix{\bfD(\A)\ar[d]^{Q_\la}\ar[rr]^-{\bfR F}&&\bfD(\B)&
\bfD(\B)\ar[d]^{Q_\r}\ar[rr]^-{(\bfR F)_\r}&&\bfD(\A)
\\
\bfK(\Inj\A)\ar[rr]^-{\hat\bfR F}&&\bfK(\Inj\B)\ar[u]_Q&
\bfK(\Inj\B)\ar[rr]^-{(\hat\bfR F)_\r}&&\bfK(\Inj\A)\ar[u]_Q\\
}$$
\end{enumerate}
\end{thm}
\begin{proof}
The composite
$$\bfK(\A)\stackrel{J_\la}\lto\bfK(\Inj\A)\stackrel{Q}\lto\bfD(\A)$$
is naturally isomorphic to the canonical functor
$\bfK(\A)\to\bfD(\A)$; see Remark~\ref{re:quotient}. Clearly, $J\comp
Q_\r$ is its right adjoint. We denote by $\bfR F$ the right derived
functor of $F$ and have $$\bfR F=Q\comp J_\la\comp \bfK(F)\comp J\comp
Q_\r.$$ Using the definition $\hat \bfR F=J_\la\comp \bfK(F)\comp J$,
we obtain $\bfR F=Q\comp \hat\bfR F\comp Q_\r$.

(1) Suppose $F$ preserves coproducts. Then $\bfK(F)$ preserves
 coproducts. It follows that $\hat\bfR F$ preserves coproducts since
 $J$ and $J_\la$ preserve coproducts. Now apply Proposition~\ref{pr:adjoints}
 to obtain a right adjoint for $\hat\bfR F$.

(2) Suppose $F$ and $\bfR F$ preserve coproducts.  Then $\bfR F$ has a
right adjoint by Proposition~\ref{pr:adjoints}. Next we show that
$$Q\comp \hat\bfR F\comp Q_\la\cong Q\comp\hat\bfR F\comp Q_\r.$$ We
have a natural transformation $Q_\la\to Q_\r$ which is induced from
the natural transformation $Q_\la\comp Q\to Q_\r\comp Q$.  Now apply
$Q\comp\hat\bfR F$ to get a natural transformation
$$\mu\colon Q\comp \hat\bfR F\comp Q_\la\lto Q\comp\hat\bfR F\comp
Q_\r.$$ It is shown in Lemma~\ref{le:nat_transf} that $Q_\la\to Q_\r$
is an isomorphism when restricted to compact objects in $\bfD(\A)$. On
the other hand, $Q\comp \hat\bfR F\comp Q_\la$ and $Q\comp\hat\bfR
F\comp Q_\r$ both preserve coproducts by our assumption on $\bfR
F$. It follows that $\mu$ is an isomorphism since $\bfD(\A)$ is
compactly generated. Clearly, $Q\comp (\hat\bfR F)_\r\comp Q_\r$ is a
right adjoint for $Q\comp \hat\bfR F\comp Q_\la$. This completes the
proof.
\end{proof}

The extended derived functor and its right adjoint admit some
alternative description. I am indebted to Bernhard Keller for
providing this remark.

\begin{rem} 
  It is possible to express $\hat\bfR F$ as the tensor functor and its
  right adjoint $(\hat\bfR F)_\r$ as the Hom functor with respect to a
  bimodule of DG categories; see \cite[6.4]{Ke}. This depends on the
  appropriate choice of DG categories $\A_0$ and $\B_0$ such that
  $\bfK(\Inj\A)\cong\bfD_\dg(\A_0)$ and
  $\bfK(\Inj\B)\cong\bfD_\dg(\B_0)$ respectively.
\end{rem}
  
Next we consider the following diagram
$$\xymatrix{ \bfK(\Inj\A)\ar[d]^{Q}\ar[rr]^-{\hat\bfR
F}&&\bfK(\Inj\B)\ar[d]^Q\\ \bfD(\A)\ar[rr]^-{\bfR F}&&\bfD(\B)}$$
and ask when it is commutative.

\begin{lem}\label{le:RF}
Keep the assumptions from Theorem~\ref{th:extensions}.  There is a
natural transformation $Q\comp\hat\bfR F\to \bfR F\comp Q$ which is an
isomorphism if and only if $F$ sends every acyclic complex of
injective objects to an acyclic complex.
\end{lem}
\begin{proof}
We apply the the localization sequence
$$\xymatrix{\bfK_\ac(\Inj\A)\ar[r]^-I&\bfK(\Inj\A)
\ar[r]^-Q&\bfD(\A)}$$ from Proposition~\ref{pr:seq}. Let $X$ be an
object in $\bfK(\Inj\A)$ and consider the triangle
$$(I\comp I_\r)X\lto X\lto(Q_\r\comp Q)X\lto\Si(I\comp I_\r)X$$ in
$\bfK(\Inj\A)$. Now apply $Q\comp\hat\bfR F$ which gives a map
$$(Q\comp\hat\bfR F)X\lto (\bfR F\comp Q)X$$ since $Q\comp \hat\bfR
F\comp Q_\r\cong \bfR F$, by Theorem~\ref{th:extensions}.
Clearly, this map is an isomorphism if and only if $Q\comp \hat\bfR
F$ annihilates $(I\comp I_\r)X$.
\end{proof}
We include a simple example which illustrates the preceding lemma.

\begin{exm} 
Let $k$ be a field and $\La=k[t]/(t^2)$. We take the functor
$$F\colon\Mod\La\lto\Mod k,\quad X\mapsto\Hom_\La(k,X),$$ and observe that the
following diagram does not commute.
$$\xymatrix{ \bfK(\Inj\La)\ar[d]^{Q}\ar[rr]^-{\hat\bfR F}&&\bfK(\Inj
k)\ar[d]_{\wr}^Q \\ \bfD(\Mod\La)\ar[rr]^-{\bfR F}&&\bfD(\Mod k) }$$
For instance, we have $QX=0$ and $({\hat\bfR F})X\neq 0$ if we take
for $X$ the acyclic complex
$$\cdots\stackrel{t}\lto\La\stackrel{t}\lto\La\stackrel{t}
\lto\La\stackrel{t}\lto \cdots$$ in $\bfK(\Inj\La)$.
\end{exm}

Now we specialize and consider as an example a morphism
$f\colon\bbX\to\bbY$ between separated noetherian schemes. Let
$f_*\colon\Qcoh\bbX\to\Qcoh\bbY$ denote the direct image functor.
Note that the right derived functor $\bfR f_*\colon
\bfD(\Qcoh\bbX)\to\bfD(\Qcoh\bbY)$ preserves coproducts
\cite[Lemma~1.4]{N1}. Thus $\bfR f_*$ and its right adjoint
Grothendieck duality functor $f^!$ extend to functors between
$\bfK(\Inj\bbX)$ and $\bfK(\Inj\bbY)$, by
Theorem~\ref{th:extensions}. This is the statement of
Theorem~\ref{th:duality} from the introduction. In fact, the situation
is in this case much nicer.  I am grateful to Amnon Neeman for
pointing out that the functor $\hat\bfR f_*$ and its right adjoint
$(\hat\bfR f_*)_\r$ make the following diagram commutative.
\begin{equation}\label{eq:RF}
\xymatrix{ \bfK(\Inj\bbX)\ar[d]^{Q}\ar[rr]^-{\hat\bfR
f_*}&&\bfK(\Inj\bbY)\ar[d]^Q&\bfK(\Inj\bbY)\ar[rr]^-{(\hat\bfR
f_*)_\r}&&\bfK(\Inj\bbX)\\ \bfD(\Qcoh\bbX)\ar[rr]^-{\bfR
f_*}&&\bfD(\Qcoh\bbY)&
\bfD(\Qcoh\bbY)\ar[u]_{Q_\r}\ar[rr]^-{f^!}&&\bfD(\Qcoh\bbX)\ar[u]_{Q_\r}
}\end{equation} This is essentially the statement of
Theorem~\ref{th:sduality} from the introduction.
The proof which is due to Neeman goes as follows.

\begin{proof}[Proof of Theorem~\ref{th:sduality}]
We need to show that both squares in (\ref{eq:RF}) commute.  Then we
use the localization sequence
$$\xymatrix{\bfS(\Qcoh\bbX)\ar[r]^-I&\bfK(\Inj\bbX)
\ar[r]^-Q&\bfD(\Qcoh\bbX)}$$ from Proposition~\ref{pr:seq} and obtain
from $\hat\bfR f_*$ and $(\hat\bfR f_*)_\r$ an adjoint pair of
functors between $\bfS(\Qcoh\bbX)$ and $\bfS(\Qcoh\bbY)$.

In order to show the commutativity of (\ref{eq:RF}), we apply
Lemma~\ref{le:RF} and need to show that $f_*$ sends an acyclic complex
$X$ of injective objects to an acyclic complex.  The question is local
in $\bbY$ and we may assume $\bbY$ affine. Cover $\bbX$ by a finite
number of affines. Then $f_*$ can be computed using the \v{C}ech
cohomology of the cover. If there are $n$ open sets in the cover, then
for any quasi-coherent sheaf $A$ we have
$$\bfR^{n+1}f_* A = 0.$$
Now take our acyclic complex $X$ of injective
sheaves on $\bbX$. Then the sequence $$0\lto X^0\lto X^1\lto
X^2\lto\cdots$$
is an injective resolution of the kernel $A$ of the
map $X^0\to X^1$.  Applying $f_*$, the sequence computes for us
$\bfR^i f_* A$, which vanishes if $i \geqslant n+1$. Thus $f_*X$ is
acyclic above degree $n$, but by shifting we conclude that it is
acyclic everywhere.  

Having shown the commutativity of the left hand square, the
commutativity of the right hand square follows, because it is
obtained by taking right adjoints.  Thus the proof is complete.
\end{proof}

Next we investigate for an exact functor $F\colon\A\to\B$ an extension
$\hat\bfL F$ of the derived functor $\bfL
F\colon\bfD(\A)\to\bfD(\B)$. For this we need some assumptions, and it
is convenient to introduce the following notation. As before, $\A$ and
$\B$ denote locally noetherian Grothendieck categories. Let
$f\colon\noeth\A\to\noeth\B$ be an additive functor. Then there is, up
to isomorphism, a unique functor $f^*\colon\A\to\B$ which extends $f$
and preserves filtered colimits. This has a right adjoint
$f_*\colon\B\to\A$ if and only if $f$ is right exact.  Note that $f$
is exact iff $f^*$ is exact iff $f_*$ sends injective objects to
injective objects.  Here is an example.

\begin{exm}\label{ex:f}
Let $f\colon\bbX\to\bbY$ be a morphism between noetherian
schemes. Then the inverse image functor
$f^*\colon\Qcoh\bbY\to\Qcoh\bbX$ sends coherent sheaves to coherent
sheaves and preserves filtered colimits. Moreover, the direct image
functor $f_*$ is a right adjoint of $f^*$. Our notation is therefore
consistent if we identify the morphism $f\colon\bbX\to\bbY$ with the
functor $\coh\bbY\to\coh\bbX$.
\end{exm}

\begin{thm}\label{th:LF}
Let $\A$ and $\B$ locally noetherian Grothendieck categories such that
$\bfD(\A)$ and $\bfD(\B)$ are compactly generated.  Let $f\colon
\noeth\A\to\noeth\B$ be an exact functor. Then $\hat\bfR f_*$ has a
left adjoint $\hat\bfL f^*$ which induces a functor $\bfS f^*$ making
the following diagram commutative.
\begin{equation}\label{eq:LF1}
\xymatrix{
&&&\bfK^c(\Inj\A)\ar@{>->}[ld] \ar[rr]^-\sim\ar'[d][dd]
&&\bfD^b(\noeth\A)\ar@{>->}[ld]\ar[dd]^-{\bfD^b(f)}\\
\bfS(\A)\ar[rr]^-I\ar[dd]^-{\bfS f^*}&&\bfK(\Inj\A)
\ar[rr]^(.7)Q\ar[dd]^-{\hat\bfL f^*}&&\bfD(\A)\ar[dd]^(.3){\bfL f^*}\\
&&&\bfK^c(\Inj\B)\ar@{>->}[ld]\ar'[r][rr]^-\sim&&\bfD^b(\noeth\B)\ar@{>->}[ld]\\
\bfS(\B)\ar[rr]^-I&&\bfK(\Inj\B) \ar[rr]^-Q&&\bfD(\B)}
\end{equation} 
If $\bfR f_*$ preserves coproducts,
then in addition the following diagram commutes.
$$\xymatrix{
&\bfD^b(\noeth\A)/\bfD^c(\A)\ar[ld]_-{F_\A}\ar'[d][dd]^-{\overline{\bfD^b(f)}}&&
\bfK^c(\Inj\A)\ar'[d][dd]\ar[ll]\ar@{>->}[ld]
&&\bfD^c(\A)\ar[dd]\ar[ll]\ar@{>->}[ld]\\
\bfS(\A)\ar[dd]^-{\bfS f^*}&&\bfK(\Inj\A)\ar[ll]_(.3){I_\la}
\ar[dd]^(.3){\hat\bfL f^*}&&\bfD(\A)\ar[ll]_(.3){Q_\la}\ar[dd]^(.3){\bfL f^*}\\
&\bfD^b(\noeth\B)/\bfD^c(\B)\ar[ld]_-{F_\B}&&
\bfK^c(\Inj\B)\ar'[l][ll]\ar@{>->}[ld]
&&\bfD^c(\B)\ar'[l][ll]\ar@{>->}[ld]\\
\bfS(\B)&&\bfK(\Inj\B)\ar[ll]_-{I_\la}&&\bfD(\B)\ar[ll]_-{Q_\la}
}$$
\end{thm}
Note that $F_\A$ and $F_\B$ induce, up to direct
factors, equivalences onto the full subcategories of compact objects
in $\bfS(\A)$ and $\bfS(\B)$ respectively. Thus $\overline{\bfD^b(f)}$
determines the functor $\bfS f^*$.

\begin{proof}
The exactness of $f$ implies the exactness of $f^*$. Thus $f_*$ sends
injective objects to injective objects and we have the following
commutative diagram.
\begin{equation}\label{eq:JKF}
\xymatrix{\bfK(\Inj\B)\;\ar@{>->}[rr]^-{J}\ar[d]^-{\hat\bfR f_*}&&
\bfK(\B)\ar[d]^-{\bfK (f_*)}\\
\bfK(\Inj\A)\;\ar@{>->}[rr]^-{J}&&\bfK(\A)}
\end{equation} 
The right adjoint $f_*$ preserves products and we have therefore a
left adjoint for $\hat\bfR f_*$, by Proposition~\ref{pr:adjoints},
which we denote by $\hat\bfL f^*$. We obtain the following diagram
$$\xymatrix{\bfK(\A)\ar[rr]^-{J_\la}\ar[d]^-{\bfK (f^*)}&&\bfK(\Inj\A)
\ar[rr]^-Q\ar[d]^-{\hat\bfL f^*}&&\bfD(\A)\ar[d]^-{\bfL f^*}\\
\bfK(\B)\ar[rr]^-{J_\la}&&\bfK(\Inj\B) \ar[rr]^-Q&&\bfD(\B)}$$ and
claim it is commutative. The left hand square commutes because it is
obtained from (\ref{eq:JKF}) by taking left adjoints. The outer square
commutes because the composite $Q\comp J_\la$ is naturally isomorphic
to the canonical functor $\bfK(\A)\to\bfD(\A)$; see
Remark~\ref{re:quotient}. We conclude the commutativity of the right
hand square, using that $J_\la\comp
J\cong\Id_{\bfK(\Inj\A)}$. Clearly, $\hat\bfL f^*$ sends acyclic
complexes to acyclic complexes and we obtain the functor $\bfS f^*$
making the diagram (\ref{eq:LF1}) commutative.  Finally, observe that
$\hat\bfL f^*$ preserve compactness because its right adjoint
$\hat\bfR f_*$ preserve coproducts; see Lemma~\ref{le:adjoint}.
Thus every square in the diagram (\ref{eq:LF1}) commutes.

Now assume in addition that $\bfR f_*$ preserves coproducts.  We use
again the fact that a left adjoint preserves compactness if the right
adjoint preserves coproducts; see Lemma~\ref{le:adjoint}.  Thus $\bfL
f^*$ and $\hat\bfL f^*$ preserve compactness. Note that $Q_\la$
identifies $\bfD(\A)$ with the localizing subcategory of
$\bfK(\Inj\A)$ which is generated by
$$\bfD^c(\A)\subseteq\bfD^b(\noeth\A)\cong
\bfK^c(\Inj\A)\subseteq\bfK(\Inj\A).$$ Of course, the same applies for
$\B$.  We obtain the following diagram
\begin{equation}\label{eq:JKF1}
\xymatrix{\bfS(\A)\ar[d]^-{S}&&\bfK(\Inj\A)
\ar[ll]_-{I_\la}\ar[d]^-{\hat\bfL
f^*}&&\bfD(\A)\ar[ll]_-{Q_\la}\ar[d]^-{\bfL f^*}\\
\bfS(\B)&&\bfK(\Inj\B) \ar[ll]_-{I_\la}&&\bfD(\B)\ar[ll]_-{Q_\la}}
\end{equation}
where the right hand square commutes. The horizontal sequences are
localization sequences by Theorem~\ref{th:products}, and $\hat\bfL
f^*$ induces a functor $S\colon\bfS(\A)\to\bfS(\B)$ making the left
hand square commutative. Moreover, we have
$$S=S\comp I_\la\comp I=I_\la\comp\hat\bfL f^*\comp I=I_\la\comp
I\comp \bfS f^*=\bfS f^*.$$ The functors $F_\A$ and $F_\B$ are both
induced by $I_\la$, and the commutativity
$$\bfS f^*\comp F_\A=F_\B\comp\overline{\bfD^b(f)}$$ is easily
checked; see Corollary~\ref{co:stable}. This completes the proof.
\end{proof}

I am grateful to the referee for pointing out possible
generalizations.

\begin{rem}
Let $f\colon \bbX\to\bbY$ be a morphism between noetherian schemes and
suppose $f^*$ is exact. Then $\bfS f^*$ is the left adjoint of $\bfS
f_*$ which appears in Theorem~\ref{th:sduality}. In fact, this theorem
suggests that parts of Theorem~\ref{th:LF} can be generalized. For
instance, the right hand square in diagram (\ref{eq:JKF1}) does not
need any assumption on the morphism $f$ because it is simply the left
adjoint of a commutative square in Theorem~\ref{th:sduality}.
\end{rem}

Next we investigate the inclusion $f\colon \bbX\to\bbY$ of an open
subscheme.  In this case, the adjoint pair of functors $f_*$ and $f^*$
between $\Qcoh\bbX$ and $\Qcoh\bbY$ restricts to an adjoint pair of
functors between $\Inj\bbX$ and $\Inj\bbY$; see \cite[VI]{G}.
Moreover, $f^*\comp f_*\cong\Id_{\Qcoh\bbX}$. Thus we can identify
$\hat\bfR f_*=f_*$ and $\hat\bfL f^*=f^*$. Note that both functors
send acyclic complexes with injective components to acyclic complexes.
This is clear for $f^*$ because it is exact, and follows for $f_*$
from Theorem~\ref{th:sduality}, or by looking at the right adjoint of
the right hand square in diagram (\ref{eq:JKF1}).  We denote for each
sheaf $A$ in $\Qcoh\bbY$ by $\Supp A$ the support of $A$ and observe
that $f^*$ annihilates $A$ if and only $\Supp A$ is contained in
$\bbY\setminus\bbX$. In fact, the natural map $A\to (f_*\comp f^*)A$
induces a split exact sequence
$$0\lto A'\lto A\lto (f_*\comp f^*)A\lto 0$$ if $A$ is injective.  In
particular, the support of $A'$ is contained in $\bbY\setminus\bbX$,
whereas the support of $(f_*\comp f^*)A$ is contained in $\bbX$.

Now fix a complex $X$ in $\bfK(\Inj\bbY)$. The {\em support} of $X$ is
by definition
$$\Supp X=\bigcup_{n\in\bbZ}\Supp X^n,$$ where $X$ is assumed to be
homotopically minimal; see Proposition~\ref{pr:min}.  We write
$X_{\bbX}=(f_*\comp f^*)X$, and the natural map $X\to X_{\bbX}$
induces an exact triangle
\begin{equation}\label{eq:supp}
X_{\bbY\setminus\bbX}\lto X\lto
X_{\bbX}\lto\Si(X_{\bbY\setminus\bbX})
\end{equation} 
in $\bfK(\Inj\bbY)$ where the support of $X_{\bbY\setminus\bbX}$ is
contained in ${\bbY\setminus\bbX}$.

\begin{lem}\label{le:supp}
  Let $\bbY$ be a seperated noetherian scheme and $f\colon
  \bbX\to\bbY$ be the inclusion of an open subscheme.  If $X$ is a
  complex in $\bfK(\Inj\bbY)$, then $f^*X=0$ if and only if the
  support of $X$ is contained in ${\bbY\setminus\bbX}$.
\end{lem}
\begin{proof} 
We have $f^*X=0$ if and only if the first map in the triangle
(\ref{eq:supp}) is an isomorphism.
\end{proof}
  
It is well known that $f^*$ induces an equivalence
$$\bfD(\Qcoh\bbY)/\bfD_{\bbY\setminus\bbX}(\Qcoh\bbY)\stackrel{\sim}\lto
\bfD(\Qcoh\bbX),$$ where $\bfD_{\bbY\setminus\bbX}(\Qcoh\bbY)$ denotes
the full subcategory of all complexes in $\bfD(\Qcoh\bbY)$ such that
the support of the cohomology is contained in $\bbY\setminus\bbX$.
We obtain an analogue for $\bfK(\Inj\bbY)$ and $\bfS(\Qcoh\bbY)$ if we define
\begin{align*}
\bfK_{\bbY\setminus\bbX}(\Inj\bbY)&=\{X\in\bfK(\Inj\bbY)\mid\Supp
X\subseteq \bbY\setminus\bbX\},\\
\bfS_{\bbY\setminus\bbX}(\Qcoh\bbY)&=\{X\in\bfS(\Qcoh\bbY)\mid\Supp
X\subseteq \bbY\setminus\bbX\}.
\end{align*}

\begin{prop}\label{pr:supp}
  Let $\bbY$ be a seperated noetherian scheme and $f\colon
  \bbX\to\bbY$ be the inclusion of an open subscheme.  Then $f^*$
  induces equivalences
\begin{eqnarray*}
&\bfK(\Inj\bbY)/\bfK_{\bbY\setminus\bbX}(\Inj\bbY)\stackrel{\sim}\lto
  \bfK(\Inj\bbX),\\
&\bfS(\Qcoh\bbY)/\bfS_{\bbY\setminus\bbX}(\Qcoh\bbY)\stackrel{\sim}\lto
  \bfS(\Qcoh\bbX).
\end{eqnarray*}
\end{prop}
\begin{proof}
  We have $f^*\comp f_*\cong\Id_{\Qcoh\bbX}$, and this carries over to
  complexes of injectives. On the other hand, we have for $X$ in
  $\bfK(\Inj\bbY)$ a natural map $X\to (f_*\comp f^*)X$ which induces
  an isomorphism in
  $\bfK(\Inj\bbY)/\bfK_{\bbY\setminus\bbX}(\Inj\bbY)$, by
  Lemma~\ref{le:supp}. This gives the first equivalence. The second
  equivalence follows from the first, because $f^*$ and $f_*$ restrict
  to functors between $\bfS(\Qcoh\bbY)$ and $\bfS(\Qcoh\bbX)$. This is
  clear for $f^*$ because it is exact. For $f_*$ this follows from
  Theorem~\ref{th:sduality}.
\end{proof}

Let us give a more elaborate formulation of Proposition~\ref{pr:supp}.
The functor $\hat\bfR f_*=f_*\colon\bfK(\Inj\bbX)\to\bfK(\Inj\bbY)$
admits a left and a right adjoint.  Therefore $\hat\bfR f_*$ induces a
recollement
$$\xymatrix{\bfK(\Inj\bbX)\ar[r]&\bfK(\Inj\bbY)
  \ar[r]\ar@<.75ex>[l]\ar@<-.75ex>[l]&\bfK_{\bbY\setminus\bbX}(\Inj\bbY).
  \ar@<.75ex>[l]\ar@<-.75ex>[l]}$$
This recollement is compatible with
the recollement
$$\xymatrix{\bfS(\Qcoh\bbY)\ar[r]&\bfK(\Inj\bbY)
\ar[r]\ar@<.75ex>[l]\ar@<-.75ex>[l]&\bfD(\Qcoh\bbY),
\ar@<.75ex>[l]\ar@<-.75ex>[l]}$$
and we obtain the following diagram.
$$\xymatrix{
\bfS(\Qcoh\bbX)\ar[d]\ar[r]&\bfK(\Inj\bbX)\ar[d]
\ar[r]\ar@<.75ex>[l]\ar@<-.75ex>[l]&\bfD(\Qcoh\bbX)\ar[d]
\ar@<.75ex>[l]\ar@<-.75ex>[l]\\
\bfS(\Qcoh\bbY)\ar[d]\ar@<.75ex>[u]\ar@<-.75ex>[u]\ar[r]&\bfK(\Inj\bbY)\ar[d]
\ar[r]\ar@<.75ex>[u]\ar@<-.75ex>[u]\ar@<.75ex>[l]\ar@<-.75ex>[l]&
\bfD(\Qcoh\bbY)\ar[d]\ar@<.75ex>[u]\ar@<-.75ex>[u]
\ar@<.75ex>[l]\ar@<-.75ex>[l]\\
\bfS_{\bbY\setminus\bbX}(\Qcoh\bbY)\ar@<.75ex>[u]\ar@<-.75ex>[u]\ar[r]&
\bfK_{\bbY\setminus\bbX}(\Inj\bbY)
\ar[r]\ar@<.75ex>[u]\ar@<-.75ex>[u]\ar@<.75ex>[l]\ar@<-.75ex>[l]&
\bfD_{\bbY\setminus\bbX}(\Qcoh\bbY)\ar@<.75ex>[u]\ar@<-.75ex>[u]
\ar@<.75ex>[l]\ar@<-.75ex>[l]
}$$
In this diagram, each row and each column is a recollement. Moreover, the diagram is
commutative if one restricts to arrows in south and east direction.
All other commutativity relations follow by taking left adjoints or
right adjoints. 

Proposition~\ref{pr:supp} tells us precisely when the inclusion of a
subscheme induces an equivalence for the stable derived category.
In \cite{O}, Orlov observed that the bounded stable derived
category of a noetherian scheme depends only on the singular
points. We extend this result to the unbounded stable derived
category, using a completely different proof.

\begin{cor}\label{co:subscheme}
  Let $\bbY$ be a seperated noetherian scheme of finite Krull
  dimension. If $f\colon \bbX\to\bbY$ denotes the inclusion of an open
  subscheme which contains all singular points of $\bbY$, then $\bfS
  f^*\colon \bfS(\Qcoh\bbY)\to\bfS(\Qcoh\bbX)$ is an equivalence.
\end{cor}
\begin{proof}
  We apply Proposition~\ref{pr:supp} and need to show that
  $\bfS_{\bbY\setminus\bbX}(\Qcoh\bbY)=0$. But this is clear from our
  assumptions on $\bbX$ and $\bbY$.
\end{proof}

Orlov's result \cite[Proposition~1.14]{O} is an immediate consequence
if one restricts the equivalence $\bfS f^*$ to compact objects; see
Theorem~\ref{th:LF}.

\begin{cor}
  Let $\bbY$ be a seperated noetherian scheme of finite Krull
  dimension. If $f\colon \bbX\to\bbY$ denotes the inclusion of an open
  subscheme which contains all singular points of $\bbY$, then $f^*$
  induces (up to direct factors) an equivalence
$$\bfD^b(\coh\bbY)/\bfD^\perf(\coh\bbY)\lto
  \bfD^b(\coh\bbX)/\bfD^\perf(\coh\bbX).$$
\end{cor}

\section{Gorenstein injective approximations and Tate cohomology}\label{se:gor}

Let $\A$ be a locally noetherian Grothendieck category and suppose
that the derived category $\bfD(\A)$ is compactly generated.  In this
section, we study the category of complete injective resolutions. We
assign functorially to each complex of injectives a complete
resolution. This yields Gorenstein injective approximations and Tate
cohomology groups for objects in $\A$.  The classical definition of
Tate cohomology is based on complete projective resolutions.  Our
approach is essentially the same, using however resolutions with
injective instead of projective components.  Another aspect in this
section is the interplay between the stable derived category
$\bfS(\A)$ and the {\em stable category} $\underline\A$ modulo
injective objects, which is obtained from $\A$ by identifying two maps
if their difference factors through some injective object.  Given
objects $A,B$ in $\A$, we write
$$\uHom_\A(A,B)=\Hom_{\underline\A}(A,B).$$ The functor
$$\bfK(\Inj\A)\lto\underline\A, \quad X\mapsto Z^0X=\Ker(X^0\to X^1)$$
provides a link between the stable categories $\bfS(\A)$ and
$\underline\A$.  In particular, we obtain an explicit description of
the stabilization functor
$$S\colon \A\xto{\can}\bfD(\A)\xto{I_\la\comp Q_\r}\bfS(\A)$$
provided
that $\A$ has some appropriate Gorenstein property. 

Most of the concepts in this section are classical, but seem to be new
in this setting and this generality. We refer to the end of this
section for historical remarks and references to the literature.

Let us start with the relevant definitions.  A complex $X$ in $\Inj\A$ is
called {\em totally acyclic} if $\Hom_\A(A,X)$ and $\Hom_\A(X,A)$ are
acyclic complexes of abelian groups for all $A$ in $\Inj\A$. We denote
by $\bfK_\tac(\Inj\A)$ the full subcategory of all totally acyclic
complexes in $\bfK(\Inj\A)$. Following \cite{EJ}, we call an object
$A$ in $\A$ {\em Gorenstein injective} if it is of the form
$Z^0X$ for some $X$ in $\bfK_\tac(\Inj\A)$.  We write
$\GInj\A$ for the full subcategory formed by all Gorenstein injective
objects.

\begin{lem}\label{le:sigma_XY}
Let $\A$ be an abelian category and let $X,Y$ be objects in
$\bfK(\Inj\A)$. Suppose $H^nX=0$ for all $n>0$ and $Y$ is totally
acyclic. Then the canonical map
$$\s_{X,Y}\colon\Hom_{\bfK(\Inj\A)}(X,Y)\lto\uHom_{\A}(Z^0X,Z^0Y)$$ is
bijective.
\end{lem}
\begin{proof}
Fix a map $\a\colon Z^0X\to Z^0Y$ in $\A$. We need to extend $\a$ to a
chain map $\bar\a\colon X\to Y$ such that $Z^0\bar\a=\a$. We use the
assumption on $X$ to extend $\a$ in non-negative degrees, and the
assumption on $Y$ allows to extend $\a$ in negative degrees. Thus
$\s_{X,Y}$ is surjective.  To show that $\s_{X,Y}$ is injective, let
$\p\colon X\to Y$ be a chain map such that $Z^0\p$ factors through
some injective object. A similar argument as before yields a chain
homotopy $X\to Y$ which shows that $\p$ is null homotopic. Thus the
proof is complete.
\end{proof}

Let us denote by $\GInj\underline\A$ the full subcategory of
$\underline\A$ formed by the objects in $\GInj\A$. Observe that
$\GInj\A$ is a Frobenius category with respect to the class of exact
sequences from $\A$. With respect to this exact structure, an object
$A$ in $\GInj\A$ is projective iff $A$ is injective iff $A$ belongs to
$\Inj\A$.  Thus the category $\GInj\underline\A$ carries a
triangulated structure. The shift takes an object $A$ to the cokernel
$\Si A$ of a monomorphism $A\to E$ into an injective object $E$.  The
exact triangles are induced from short exact sequences in $\A$.

\begin{prop}\label{pr:gproperty}
Let $\A$ be an abelian category. Then the functor
$$\bfK_\tac(\Inj\A)\lto\GInj\underline\A, \quad X\mapsto Z^0X,$$
is an equivalence of triangulated categories.
\end{prop}
\begin{proof}
We need to show that the functor is fully faithful and surjective on
isomorphism classes of objects. The last property is clear from the
definition of $\GInj\A$. The functor is fully faithful by
Lemma~\ref{le:sigma_XY}. Finally, observe that an exact triangle of
complexes comes, up to isomorphism, from a sequence of
complexes which is split exact in each degree. Thus we obtain an exact
sequence in $\A$ and an exact triangle in $\underline\A$ if we apply $Z^0$.
\end{proof}

The following lemma is crucial because it provides the existence of
complete injective resolutions. Let us write
$$G\colon\bfK_\tac(\Inj\A)\lto\bfK(\Inj\A)$$ for the inclusion
functor.

\begin{lem}\label{le:hac}
Let $\A$ be a locally noetherian Grothendieck category and suppose
that $\bfD(\A)$ is compactly generated. Then the inclusion
$G\colon\bfK_\tac(\Inj\A)\to\bfK(\Inj\A)$ has a left adjoint
$$G_\la\colon\bfK(\Inj\A)\lto\bfK_\tac(\Inj\A).$$
\end{lem}
\begin{proof}
  The inclusion $I\colon\bfK_\ac(\Inj\A)\to\bfK(\Inj\A)$ has a left
  adjoint $I_\la$ by Theorem~\ref{th:products}. Thus it is sufficient
  to show that the inclusion $\bfK_\tac(\Inj\A)\to\bfK_\ac(\Inj\A)$
  has a left adjoint.  Let us denote by $E$ the coproduct of a
  representative set of all indecomposable injective objects in
  $\A$. By definition, $\bfK_\tac(\Inj\A)$ consists of all objects $X$
  in $\bfK_\ac(\Inj\A)$ such that
$$\Hom_{\bfK(\Inj\A)}(\Si^nE,X)\cong\Hom_{\bfK_\ac(\Inj\A)}(I_\la(\Si^nE),X)$$
vanishes for all $n\in\bbZ$. The category $\bfK_\ac(\Inj\A)$ is
compactly generated by Corollary~\ref{co:stable}, and we can apply
Proposition~\ref{pr:adjoints2} to obtain a left adjoint for the
inclusion $\bfK_\tac(\Inj\A)\to\bfK_\ac(\Inj\A)$.
\end{proof}

Given an object $A$ in $\A$ with injective resolution $iA$, we call
the natural map $$iA\lto G_\la iA$$ a {\em complete injective
resolution} of $A$. If we apply the functor $Z^0$ to this map, we
obtain a Gorenstein injective approximation of $A$.

\begin{thm}\label{th:GInj}
Let $\A$ be a locally noetherian Grothendieck category and suppose
that $\bfD(\A)$ is compactly generated. Then the inclusion
$\GInj\underline\A\to\underline\A$ has a left adjoint
$$T\colon \underline\A\lto\GInj\underline\A.$$ Thus we have for each
object $A$ in $\A$ a natural map $A\to TA$ which induces a bijection
$$\uHom_\A(TA,B)\stackrel{\sim}\lto\uHom_\A(A,B)\quad\textrm{for
all}\quad B\in\GInj\A.$$
\end{thm}
\begin{proof}
Fix  an object $A$ in $\A$ and choose an
injective resolution $iA$. We put
$$TA=Z^0(G_\la iA),$$ and this induces a functor
$T\colon\underline\A\to\GInj\underline\A$. Let $B$ in $\GInj\A$ and fix a
totally acyclic complex $tB$ such that $Z^0tB=B$.  The natural map $iA\to
G_\la iA$ induces a map $A\to TA$ in $\underline\A$ which makes the
following square commutative.
$$
\xymatrix{\Hom_{\bfK(\Inj\A)}(G_\la iA,tB)\ar[d]^{Z^0}\ar[r]^-\sim&
\Hom_{\bfK(\Inj\A)}(iA,tB)\ar[d]^{Z^0}\\
\uHom_\A(TA,B)\ar[r]&\uHom_\A(A,B)}
$$ The vertical maps are bijective by Lemma~\ref{le:sigma_XY}, and we
conclude that $T$ is a left adjoint for the inclusion
$\GInj\underline\A\to\underline\A$. 
\end{proof}

Next we use complete injective resolutions to define Tate cohomology
groups for objects in $\A$.

\begin{defn} 
Given objects $A,B$ in $\A$ and $n\in\bbZ$, the {\em Tate
cohomology group} is
\begin{equation*}\label{eq:tate}
\tExt^n_\A(A,B)=H^n\Hom_\A(A,G_\la iB)
\end{equation*} 
\end{defn}

\begin{rem} 
The correct term for this cohomology theory would be `injective Tate
cohomology' in order to distinguish it from the usual `projective Tate
cohomology' which is defined via complete projective resolutions. For
simplicity, we drop the extra adjective `injective'. Note that
confusion is not possible because we do not consider projective Tate
cohomology in this paper.
\end{rem}

Tate cohomology is natural in both arguments because the formation of
complete injective resolutions is functorial. In addition, we have a
{\em comparison map}
$$\Ext^n_\A(A,B)\lto\tExt^n_\A(A,B),$$ which is
induced by the map $iB\to G_\la iB$. There is an alternative
description of Tate cohomology which is based on the left adjoint
$T\colon \underline\A\to\GInj\underline\A$.

\begin{prop}
Given objects $A,B$ in $\A$ and $n\in\bbZ$, there is a natural isomorphism
$$\tExt^n_\A(A,B)\cong\uHom_\A(A,\Si^n(TB)).$$
\end{prop}
\begin{proof}
Using Lemmas~\ref{le:com} and \ref{le:sigma_XY}, we have the following
sequence of isomorphisms
\begin{eqnarray*}
H^n\Hom_\A(A,G_\la iB)&\cong&\Hom_{\bfK(\A)}(A,\Si^n(G_\la iB))\\
&\cong&\Hom_{\bfK(\A)}(iA,\Si^n(G_\la iB))\\
&\cong&\uHom_\A(A,\Si^n(TB)).
\end{eqnarray*}
\end{proof}

We have a more conceptual
definition of Tate cohomology for objects in $\bfD(\A)$ which uses the
composite
$$U\colon\bfD(\A)\xto{G_\la\comp Q_\r}\bfK_\tac(\Inj\A).$$
Thus we define for objects $X,Y$ in $\bfD(\A)$ and $n\in\bbZ$
$$\tExt^n_\A(X,Y)=\Hom_{\bfK(\A)}(UX,\Si^n(UY)).$$ Note that this
definition is consistent with the original definition of Tate
cohomology if we take objects in $\A$ and view them as complexes
concentrated in degree $0$. This follows from the fact that $Q_\r X$
is nothing but an injective resolution of $X^0$ when $X$ is
concentrated in degree $0$. From now on, we will use one of the
alternative descriptions of Tate cohomology whenever this is
convenient.

Next we show that each exact sequence in $\A$ induces a long exact
sequence in Tate cohomology. This is based on the following simple lemma.

\begin{lem}\label{le:g1}
The left adjoint $T\colon \underline\A\to\GInj\underline\A$ has the
following properties.
\begin{enumerate}
\item An exact sequence $0\to A'\to A \to A''\to 0$ in $\A$ induces an
exact triangle
$$TA'\lto TA\lto TA''\lto\Si(TA')\quad\textrm{in}\quad \GInj\underline\A.$$
\item Let $A,B$ be in $\A$ and $n\in\bbZ$. The natural map $A\to TA$
induces an isomorphism
$$\tExt^n_\A(TA,B)\cong\tExt^n_\A(A,B).$$
\end{enumerate}
\end{lem}
\begin{proof} 
(1) We have an exact triangle $A'\to A\to A''\to\Si (A')$ in
$\bfD(\A)$. Now use that the exact functor $Z^0\comp G_\la\comp Q_\r$
computes $T$.

(2) The adjointness property of $T$ implies
    $\uHom_\A(TA,TB)\cong\uHom_\A(A,TB)$.
\end{proof}

\begin{prop}\label{pr:g2}
Let $0\to B'\to B \to B''\to 0$ be an exact sequence in $\A$. Then we
have for $A$ and $C$ in $\A$ the following long exact sequences.
\begin{multline*}
\cdots\lto\tExt^n_\A(A,B')\lto\tExt^n_\A(A,B)\lto\tExt^n_\A(A,B'')\lto\\
\lto\tExt^{n+1}_\A(A,B')\lto\tExt^{n+1}_\A(A,B)\lto\tExt^{n+1}_\A(A,B'')
\lto\cdots
\end{multline*}
\begin{multline*}
\cdots\lto\tExt^n_\A(B'',C)\lto\tExt^n_\A(B,C)\lto\tExt^n_\A(B',C)\lto\\
\lto\tExt^{n+1}_\A(B'',C)\lto\tExt^{n+1}_\A(B,C)\lto\tExt^{n+1}_\A(B',C)
\lto\cdots
\end{multline*}
\end{prop}
\begin{proof}
We apply Lemma~\ref{le:g1} and use the fact that $\uHom_\A(TA,-)$ and
$\uHom_\A(-,TC)$ are cohomological functors.
\end{proof}

We compute Tate cohomology  for Gorenstein injective objects.

\begin{prop}\label{pr:g3}
Let $A,B$ be objects in $\A$ and suppose $B$ is Gorenstein
injective. Then the comparison map
$$\Ext^n_\A(A,B)\lto\tExt^n_\A(A,B)$$ is an isomorphism for $n>0$ and
induces an isomorphism
$$\uHom_\A(A,B)\stackrel{\sim}\lto\tExt^0_\A(A,B)\quad\textrm{for}\quad n=0.$$
\end{prop}
\begin{proof} 
Our assumption implies $TB= B$. The case $n=0$ is clear. For $n=1$,
choose an exact sequence $0\to B\to E\to \Si B\to 0$ with $E$
injective and apply $\Hom_\A(A,-)$. The cokernel of
$$\Hom_\A(A,E)\lto\Hom_\A(A,\Si B)$$ is isomorphic $\Ext^1_\A(A,B)$;
it is isomorphic to $\uHom_\A(A,\Si B)$ since $B$ is Gorenstein
injective.  For $n>1$, use dimension shift.
\end{proof}

Next we describe those objects $A$ in $\A$ such that $\tExt^*_\A(A,-)$
vanishes. For instance, Tate cohomology vanishes for all objects having
finite projective or finite injective dimension.

\begin{prop}\label{pr:g4}
For an object $A$ in $\A$, the following are equivalent.
\begin{enumerate}
\item $\tExt^*_\A(A,-)=0$.
\item $\tExt^0_\A(A,A)=0$.
\item $\tExt^*_\A(-,A)=0$.
\item $\Ext^1_\A(A,B)=0$ for all $B\in\GInj\A$.
\item $\uHom_\A(A,B)=0$ for all $B\in\GInj\A$.
\end{enumerate}
\end{prop}
\begin{proof} 
Use the isomorphism $\tExt^*_\A(-,-)\cong\uHom^*_\A(T-,T-)$ and
the fact that $\tExt^1_\A(-,B)$ can be computed via $\Ext^1_\A(-,TB)$.
\end{proof}

The following result formulates our analogue of the maximal
Cohen-Macaulay approximation in the sense of Auslander and Buchweitz
\cite{AB}. Note that a Gorenstein injective object is the `dual' of a
maximal Cohen-Macaulay object which one defines in a category with
enough projectives. Let
$$\X=\{A\in\A\mid\Ext^1_\A(A,B)=0 \textrm{ for all }
B\in\GInj\A\}\quad\textrm{and}\quad\Y=\GInj\A.$$

\begin{thm}\label{th:cotorsion}
  Let $\A$ be a locally noetherian Grothendieck category and suppose
  that $\bfD(\A)$ is compactly generated. 
\begin{enumerate}
\item Every object $A$ in
  $\A$ fits into exact sequences
\begin{equation*}\label{eq:app}
0\to Y_A\to X_A\to A\to 0\quad\textrm{and}\quad0\to A\to Y^A\to X^A\to
0
\end{equation*} 
in $\A$ with $X_A,X^A$ in $\X$ and $Y_A,Y^A$ in
$\Y$.
\item The map $A\mapsto X_A$ induces a right adjoint for the inclusion
  $\underline\X\to\underline\A$.
\item The map $A\mapsto Y^A$ induces a left adjoint for the inclusion
  $\underline\Y\to\underline\A$.
\item $\X\cap\Y=\Inj\A$.
\end{enumerate}
\end{thm} 

Note that this is essentially the statement of Theorem~\ref{th:tate}
from the introduction, since $\X$ is precisely the subcategory of
objects $A$ in $\A$ such that the Tate cohomology functor
$\tExt_\A^*(A,-)$ vanishes.

\begin{proof}
We use the basic properties of Tate cohomology.

(1) Fix an object $A$ in $\A$ and a complete injective
resolution $$iA\lto G_\la iA=yA.$$ We complete this map to an exact
triangle
\begin{equation}\label{eq:app2}
iA\lto yA\lto xA \lto\Si(iA).
\end{equation}
in $\bfK(\Inj\A)$ and have therefore a sequence $0\to iA\to yA\to
xA\to 0$ of complexes which is split exact in each degree. Applying
$Z^0\colon\bfK(\Inj\A)\to\underline\A$ produces an exact
sequence
$$0\lto A\stackrel{\a}\lto Y^A\stackrel{\b}\lto X^A\lto 0$$ in
$\A$. Clearly, $Y^A$ belongs to $\Y$. On the other hand,
$\tExt^*_\A(\a,-)$ is an isomorphism. Thus $\tExt_\A^*(X^A,-)$
vanishes and $X^A$ belongs to $\X$.  The second sequence ending in $A$
is obtained by rotating the triangle (\ref{eq:app2}).

(2) We consider the exact sequence
$$0\lto Y_A\stackrel{\mu}\lto X_A\stackrel{\nu}\lto A\lto 0$$ and need
to show that $\uHom_\A(X,\nu)$ is bijective for all $X$ in $\X$.  To
see this, let $\p\colon X\to A$ be a map with $X$ in $\X$. The map
$\p$ factors through $\nu$ since $\Ext^1_\A(X,Y_A)=0$.  Therefore
$\uHom_\A(X,\nu)$ is surjective. To show that $\uHom_\A(X,\nu)$ is
injective, let $\psi\colon X\to X_A$ be a map such that $\nu\comp\psi$
has a factorization $X\stackrel{\p'}\to E\stackrel{\p''}\to A$ with
injective $E$. We obtain a factorization $\p''=\nu\comp\chi$, since
$E$ belongs to $\X$. We have $\nu\comp(\psi-\chi\comp\p')=0$, and
$\psi-\chi\comp\p'$ needs to factor through $\mu$. Therefore
$\psi-\chi\comp\p'$ factors through some injective object, since
$\uHom_\A(X,Y_A)=0$. We conclude that $\psi$ factors through an
injective object. Thus the map $\uHom_\A(X,\nu)$ is bijective.

(3) We consider the exact sequence
$$0\lto A\stackrel{\a}\lto Y^A\stackrel{\b}\lto X^A\lto 0$$ and need
to show that $\uHom_\A(\a,Y)$ is bijective for all $Y$ in $\Y$.  But
this is clear from the long exact sequence for Tate cohomology, since
$\uHom_\A(-,Y)\cong\tExt_\A^0(-,Y)$ and $\tExt_\A^*(X^A,-)$ vanishes.

(4) Clearly, $\Inj\A$ is contained in $\X\cap\Y$. Now let $A$ be in
$\X\cap\Y$. Thus
$$\uHom_\A(A,A)\cong\tExt_\A^0(A,A)=0.$$ Therefore the identity map
$A\to A$ factors through an injective object. We conclude that $A$ is
injective.
\end{proof}

Let us comment on the interplay between the stable category
$\underline\A$ and the stable derived category $\bfS(\A)$. We have
already seen that the definition of Tate cohomology is possible in
both settings. It is more elementary in $\underline\A$, but more
conceptual using the category of complete injective resolutions
$\bfK_\tac(\Inj\A)$ which is a subcategory of $\bfS(\A)$.  The same
phenomenon appears when one studies Gorenstein injective
approximations. The proof of Theorem~\ref{th:cotorsion} we have given
uses the category of complexes $\bfK(\Inj\A)$.  There is an
alternative proof which avoids complexes and uses instead the left
adjoint $T\colon\underline\A\to\GInj\underline\A$.

Gorenstein rings and schemes play an important role in applications
and have a number of interesting homological properties.  It is
therefore important to formulate a Gorenstein property for a locally
noetherian Grothendieck category $\A$. Let us denote by $\Si^{\infty}\A$
the full subcategory of objects $A$ in $\A$ which fit into an exact
sequence
$$\cdots \lto E_2\lto E_1\lto E_0\lto A\lto 0$$ with $E_n$ injective
for all $n$. We say that $\A$ has the {\em injective Gorenstein
property} if the equivalent conditions in the following proposition
are satisfied. This property has been studied by Beligiannis in \cite{B}.

\begin{prop}\label{pr:gorinj}
Let $\A$ be a locally noetherian Grothendieck category and suppose
that $\bfD(\A)$ is compactly generated. Then the following are
equivalent.
\begin{enumerate}
\item $\Ext^1_\A(A,B)=0$ for all $A\in\Inj\A$ and $B\in\Si^{\infty}\A$.
\item Every acyclic complex in $\Inj\A$ is totally acyclic.
\item $\GInj\A=\Si^{\infty}\A$.
\item $S\colon \A\xto{\can}\bfD(\A)\xto{I_\la\comp Q_\r}\bfS(\A)$
annihilates all injective objects. 
\item $S$ induces an equivalence
$\GInj\underline\A\to\bfS(\A)$.
\end{enumerate}
\end{prop}
\begin{proof}
The conditions (1) - (3) are pairwise equivalent. This follows from the formula
$$\Ext^n_\A(A,Z^{0}Y)\cong\Hom_{\bfK(\A)}(A,\Si^nY)\cong
H^n\Hom_\A(A,Y)$$ where $A$ is any object in $\A$ and $Y$ is an
acyclic complex in $\Inj\A$. The first isomorphism is valid for all
$n\geqslant 1$, and the second for all $n\in\bbZ$. 

Now observe that $I_\la$ annihilates precisely those objects $X$ in
$\bfK(\Inj\A)$ such that $\Hom_{\bfK(\Inj\A)}(X,Y)=0$ for every acyclic
complex $Y$ in $\Inj\A$.  On the other hand, $\A\to\bfD(\A)$ and
$Q_\r$ are faithful. Thus (1) - (3) are equivalent to (4). Also, (5)
implies (4). So, it remains to show that (1) - (4) imply (5).

Suppose (2) and (4) hold. We have already seen in
Proposition~\ref{pr:gproperty} that
$$\bfS(\A)\lto\GInj\underline\A,\quad X\mapsto Z^0X,$$ is an
equivalence, since every acyclic complex is totally acyclic.  On the
other hand, $S$ annihilates all injective objects and induces
therefore a functor $\underline\A\to\bfS(\A)$.  The composite with
$Z^0\colon\bfS(\A)\to\GInj\underline\A$ is precisely the right adjoint
of the inclusion $\GInj\underline\A\to\underline\A$ constructed in
Theorem~\ref{th:GInj}. Thus $(Z^0\comp S)A\cong A$ for all $A$ in
$\GInj\A$.
\end{proof}

\begin{cor}
Let $\A$ be a locally noetherian Grothendieck category and suppose
that $\bfD(\A)$ is compactly generated. If $\A$ has the injective
Gorenstein property, then the composite
$$\GInj\A\xto{\inc} \A\xto{\can}\bfD(\A)\xto{I_\la\comp
Q_\r}\bfS(\A)\xto{Z^0}\GInj\underline\A$$ is naturally isomorphic to
the canonical projection $\GInj\A\to\GInj\underline\A$.
\end{cor}

We are now in the position that we can describe the stabilization
functor $S\colon\A\to\bfS(\A)$, provided that $\A$ has the injective
Gorenstein property. We use the left adjoint
$T\colon\underline\A\to\GInj\underline\A$ of the inclusion
$\GInj\underline\A\to\underline\A$. For $A$ in $\A$, choose any
acyclic complex $X$ of injective objects such that $Z^0X\cong TA$.
Then $SA\cong X$.

\begin{exm}
  Let $\La$ be a ring and suppose $\La$ is {\em Gorenstein}, that is,
  $\La$ is two-sided noetherian and $\La$ has finite injective
  dimension as left and right $\La$-module. In this case, the category
  $\Mod\La$ has the injective Gorenstein property. This follows from
  the fact that every injective $\La$-module has finite projective
  dimension if $\La$ is Gorenstein; see Example~\ref{ex:pd}. Given a
  $\La$-module $A$, Tate cohomology $\tExt_\La^*(A,-)$ vanishes iff
  $A$ has finite injective dimension iff $A$ has finite projective
  dimension. Note that for Gorenstein rings, the classical Tate
  cohomology defined via complete projective resolutions coincides
  with our Tate cohomology, which is defined via complete injective
  resolutions; see \cite{BR}.
\end{exm} 

\begin{exm}\label{ex:selfinj}
Let $\La$ be a ring and suppose that projective and injective
$\La$-modules coincide. Then every $\La$-module is Gorenstein
injective. In particular, $\bfS(\Mod\La)$ is equivalent to the stable
category $\underline\A$ of $\A=\Mod\La$. Given a $\La$-module
$A$, there is an exact triangle
$$pA\lto iA\lto tA\lto \Si(pA)$$ in $\bfK(\Inj\A)$ where $pA$ denotes
a projective, $iA$ an injective, and $tA$ a Tate resolution of
$A$. This triangle is isomorphic to the canonical triangle
$$(Q_\la\comp Q)\bar A\lto \bar A\lto (I\comp I_\la)\bar A\lto
\Si(Q_\la\comp Q)\bar A$$ where $\bar A=Q_\r A$. 
\end{exm}

\begin{exm}
  Let $\bbX$ be a noetherian scheme and suppose that every injective
  object $E$ in $\Qcoh\bbX$ admits a finite resolution
  $$0\lto L_r\lto\cdots\lto L_{2}\lto L_1\lto L_0\lto E\lto 0$$
  with
  $L_n$ locally free for each $n$. Then the category $\Qcoh\bbX$ has
  the injective Gorenstein property.
\end{exm}

\subsection*{Historical remarks}
Gorenstein injective approximations and Tate cohomology have a long
history. Auslander and Bridger \cite{AuBr} introduce the stable module
category and assign to each module a G-dimension. Over Gorenstein
rings, the modules of G-dimension $0$ are precisely the maximal
Cohen-Macaulay or Gorenstein projective modules. Auslander and
Buchweitz establish maximal Cohen-Macaulay approximations in
\cite{AB}, and there is an alternative unpublished approach by
Buchweitz \cite{Bu} which involves the derived category.  Enochs and
his collaborators drop finiteness conditions on modules and prove the
existence of Gorenstein projective and Gorenstein injective
approximations for arbitrary modules, for instance over Gorenstein
rings \cite{EJ1}.  Further generalizations can be found in work of
Beligiannis \cite{B}.  J{\o}rgensen \cite{J} constructs Gorenstein
projective approximations for artin algebras via Bousfield
localization, using the category of complete projective resolutions.
Papers of Hovey \cite{H} and Beligiannis and Reiten \cite{BR} employ the
formalism of model category structures and cotorsion pairs. 

The exposition of Buchweitz \cite{Bu} discusses the close connection
between maximal Cohen-Macaulay approximations and Tate cohomology over
Gorenstein rings. For more general settings, we refer to the work of
Beligiannis and Reiten \cite{BR}. A paper of Mislin explains Tate
cohomology via satellites \cite{M}. A comparison of Tate cohomology
via projectives and injectives is carried out in work of Nucinkis
\cite{Nu}. Another exposition of Tate cohomology over noetherian rings
can be found in a paper of Avramov and Martsinkovsky \cite{AM}.

\section{Tensor products in modular representation theory}\label{se:kG}

Let $G$ be a finite group and $k$ be a field. The stable module
category $\uMod kG$ of the group algebra $kG$ plays an important role
in modular representation theory. In this section, we show that this
category is equivalent to the stable derived category $\bfS(\Mod kG)$
and study its tensor product. 

It is convenient to work in a slightly more general setting. Thus we
fix a finite dimensional cocommutative Hopf algebra $\La$ over a field
$k$ and consider the module category $\A=\Mod\La$. Note that
projective and injective modules over $\La$ coincide. The tensor
product $\otimes_k$ over $k$ induces a tensor product on $\A$ which
extends to a tensor product on $\bfK(\A)$. Similarly, $\Hom_k(-,-)$
induces products on $\A$ and $\bfK(\A)$. Note that we have a natural
isomorphism
\begin{equation}\label{eq:adj}
\Hom_{\bfK(\A)}(X\otimes_kY,Z)\cong\Hom_{\bfK(\A)}(X,\Hom_k(Y,Z))
\end{equation}
for all $X,Y,Z$ in $\bfK(\A)$. The subcategories $\bfK(\Inj\A)$ and
$\bfS(\A)$ inherit tensor products from $\bfK(\A)$ because of the
following elementary fact.

\begin{lem} 
The subcategories $\bfK(\Inj\A)$ and $\bfS(\A)$ are tensor ideals in
$\bfK(\A)$. More precisely,
\begin{enumerate} 
\item $X\in\bfK(\Inj\A)$ and $Y\in\bfK(\A)$ imply $X\otimes_k
Y\in\bfK(\Inj\A)$;
\item $X\in\bfS(\A)$ and $Y\in\bfK(\A)$ imply $X\otimes_k
Y\in\bfS(\A)$.
\end{enumerate}
\end{lem}

Now consider $k$ as a $\La$-module and view it as a complex
concentrated in degree zero; it is the unit of the tensor product in
$\bfK(\A)$. This complex fits into exact triangles
$$ak\lto k\lto ik\lto\Si(ak)\quad\textrm{and}\quad pk\lto ik\lto
tk\lto\Si(pk)$$ in $\bfK(\A)$, where $ik$ denotes an injective
resolution and $pk$ a projective resolution of $k$ in $\A$. We
consider the canonical maps $k\to ik$ and $pk\to ik$. Thus $ak$ and
$tk$ are acyclic complexes. In fact, $tk$ is a {\em Tate resolution}
of $k$ which is obtained by splicing together $pk$ and $ik$.

We have seen in previous sections that the inclusions
$$\bfS(\A)\lto\bfK(\Inj\A)\quad\textrm{and}\quad
\bfK(\Inj\A)\lto\bfK(\A)$$ have left adjoints. Next we provide
explicit descriptions of these adjoints. In particular, we see that
they preserve the tensor product $\otimes_k$.  Note that Hovey,
Palmieri, and Strickland pointed out the relevance of these categories
in their work on axiomatic stable homotopy theory \cite{HPS}; we refer
to this work for further details and applications.

\begin{thm}
Let $\La$ be a finite dimensional cocommutative Hopf algebra over a
field $k$, and let $\A=\Mod\La$.
\begin{enumerate}
\item The functor
$$-\otimes_k ik\colon\bfK(\A)\lto\bfK(\Inj\A)$$ is a left adjoint for
the inclusion $\bfK(\Inj\A)\to\bfK(\A)$.
\item The functor
$$-\otimes_k tk\colon\bfK(\Inj\A)\lto\bfS(\A)$$ is a left adjoint for
the inclusion $\bfS(\A)\to\bfK(\Inj\A)$.
\item The functor
$$\A\xto{\can}\bfK(\A)\xto{-\otimes tk} \bfS(\A)$$
induces an equivalence $\underline\A\to\bfS(\A)$ with quasi-inverse
$Z^0\colon\bfS(\A)\to \underline\A$.
\end{enumerate}
\end{thm}
\begin{proof} 
(1) Fix an object $X$ in $\bfK(\A)$. The map
$X\cong X\otimes_k k\to X\otimes_k ik$ induces for all $Y$ in $\bfK(\Inj\A)$
an isomorphism
$$\Hom_{\bfK(\A)}(X\otimes_k ik,Y)\lto\Hom_{\bfK(\A)}(X,Y).$$ This
follows from Lemma~\ref{le:com}, using in addition the formula
(\ref{eq:adj}).

(2)  Fix an object $X$ in $\bfK(\Inj\A)$. The map
$X\cong X\otimes_k ik\to X\otimes_k tk$ induces for all $Y$ in $\bfS(\A)$
an isomorphism
$$\Hom_{\bfK(\A)}(X\otimes_k tk,Y)\lto\Hom_{\bfK(\A)}(X,Y).$$ 
To see this, consider the exact triangle
$$X\otimes_k pk\lto X\otimes_k ik\lto X\otimes_k tk\lto\Si(X\otimes_k
pk).$$ Now use that $$\Hom_{\bfK(\A)}(X\otimes_k pk,Y)\cong
\Hom_{\bfK(\A)}(pk,\Hom_k(X,Y))=0$$ since $\Hom_k(X,Y)$ is acyclic.

(3) We apply Proposition~\ref{pr:gorinj}. First observe that every
object in $\A$ is Gorenstein injective.  The functor
$$\A\xto{\can}\bfK(\A)\xto{-\otimes tk} \bfS(\A)$$ is naturally
isomorphic to the stabilization functor $S\colon\A\to\bfS(\A)$. This
follows from the fact that $A\otimes_kik$ is an injective resolution
in $\A$ for each object $A$. Thus $$SA\cong (A\otimes_k ik)\otimes_k
tk\cong A\otimes_k tk.$$ In Proposition~\ref{pr:gorinj}, it is shown
that $S$ induces an equivalence $\underline\A\to\bfS(\A)$, with
quasi-inverse $Z^0$.
\end{proof}

\begin{rem} 
The unit of the product in $\bfK(\Inj\A)$ is $ik$, and its graded
endomorphism ring is the cohomology ring $H^*(\La,k)$. The unit of the
product in $\bfS(\A)$ is $tk$, and its graded endomorphism ring is
the Tate cohomology ring $\widehat{H}^*(\La,k)$.
\end{rem}

\begin{appendix}
\section{The DG category of noetherian objects}\label{ap:dg}

Let $\A$ be a locally noetherian Grothendieck category.  We give an
alternative description of the homotopy category $\bfK(\Inj\A)$ as the
derived category of some DG category. Here, we follow closely Keller's
exposition in \cite{Ke}.

Let $\C$ be a small DG category.  We recall the definition of the
derived category $\bfD_\dg(\C)$ of $\C$.  The category
$\bfC_\dg(\C)$ of {\em cochain complexes} has by definition as objects
all DG $\C$-modules. A map in $\bfC_\dg(\C)$ is a map of DG
$\C$-modules which is homogeneous of degree zero and commutes with
the differential. The {\em homotopy category} $\bfK_\dg(\C)$ is
obtained from $\bfC_\dg(\C)$ by identifying homotopy equivalent maps,
where $f,g\colon X\to Y$ are homotopy equivalent if there exists a map
$s\colon X\to Y$ of graded modules which is homogeneous of degree $-1$
and satisfies
$$(f-g)^n=s^{n+1}\comp d+d\comp s^n\quad\textrm{for all}\quad
n\in\mathbb Z.$$ Finally, the {\em derived category} of $\C$ is obtained from
$\bfK_\dg(\C)$ as the localization
$$\bfD_\dg(\C)=\bfK_\dg(\C)[Q^{-1}]$$ with respect to the class $Q$ of
all maps $f$ which induce an isomorphism $H^*f$.

Given two cochain complexes $X$ and $Y$ in $\A$, we define the
cochain complex $\HOM_\A(X,Y)$. The $n$th component is
$$\prod_{p\in\mathbb Z}\Hom_\A(X^p,Y^{n+p})$$ 
and the  differential is given by 
$$d(f^p)=d\comp f^p-(-1)^nf^{p+1}\comp d.$$

Now fix a class $\C$ of objects in $\A$. We obtain a DG category
$\bar\C$ by taking as objects for each $A$ in $\C$ an injective
resolution $\bar A$, and as maps
$$\Hom_{\bar\C}(\bar A,\bar B)=\HOM_\A(\bar A,\bar B).$$

\begin{prop}\label{pr:htp}
Let $\A$ be a locally noetherian Grothendieck category, and let $\C$
be a class of noetherian objects which generate $\bfD^b(\noeth\A)$,
that is, there is no proper thick subcategory containing $\C$.
Then the functor
$$\bfK(\Inj\A)\lto\bfD_\dg(\bar\C),\quad X\mapsto \HOM_\A(-,X)|_{\bar\C},$$
is an equivalence of triangulated categories.
\end{prop}
\begin{proof}
The functor is exact. To see that it preserves coproducts, fix an
object $A$ in $\C$ and a family of objects $X_i$ in
$\bfK(\Inj\A)$. Then we have for every $n\in\bbZ$
\begin{eqnarray*}
H^n\coprod_i\HOM_\A(\bar A,X_i)&\cong& \coprod_i H^n\HOM_\A(\bar
A,X_i) \cong \coprod_i\Hom_{\bfK(\Inj\A)}(\Si^{-n}\bar A,X_i)\\
&\cong& \Hom_{\bfK(\Inj\A)}(\Si^{-n}\bar A,\coprod_iX_i)\cong
H^n\HOM_\A(\bar A,\coprod_iX_i)
\end{eqnarray*}
since $\bar A$ is compact in $\bfK(\Inj\A)$ by Lemma~\ref{le:com}.
Thus the canonical map
$$\coprod_i\HOM_\A(-,X_i)|_{\bar\C}\lto\HOM_\A(-,\coprod_iX_i)|_{\bar\C}$$
is an isomorphism. Furthermore, the functor induces for objects $A$
and $B$ in $\C$ bijections
$$\Hom_{\bfK(\Inj\A)}(\bar A,\Si^n\bar B)\cong H^n\HOM_\A(\bar A,\bar
B)\cong H^n\Hom_{\bar\C}(\bar A,\bar B)\cong \Hom_{\bfD_\dg(\bar\C)}(\bar
A^\wedge,\Si^n\bar B^\wedge),$$ where $\bar A^\wedge$ denotes the free
module $\Hom_{\bar\C}(-,\bar A)$. Using infinite d\'evissage, we
conclude that the functor is fully faithful since $\C$ generates
$\bfK(\Inj\A)$. The functor is, up to isomorphism, surjective on
objects since the image contains the free $\bar\C$-modules which
generate $\bfD_\dg(\bar\C)$.
\end{proof}

\begin{cor}
Viewing $\noeth\A$ as DG category, we have an equivalence
$$\bfK(\Inj\A)\stackrel{\sim}\lto\bfD_\dg(\noeth\A).$$
\end{cor}

We remark that the proof of Proposition~\ref{pr:htp} works for any
homotopy category. To be precise, let $\X$ be an additive category
with arbitrary coproducts and let $\C$ be a set of objects in
$\bfK(\X)$ which are compact (when viewed as objects in the localizing
subcategory generated by $\C$). Define $\bar\C$ as before by
$$\Hom_{\bar\C}(A,B)=\HOM_\X(A,B)$$ for $A$ and $B$ in $\C$.  Then the
functor
$$\bfK(\X)\lto\bfD_\dg(\bar\C),\quad X\mapsto
\HOM_\X(-,X)|_{\bar\C},$$ induces an equivalence between the
localizing subcategory which is generated by $\C$, and
$\bfD_\dg(\bar\C)$.

\section{Homotopically minimal complexes}\label{ap:min}

A complex $X$ in some additive category is called {\em homotopically
  minimal}, if every map $\p\colon X\to X$ of complexes is an
isomorphism provided there is a map $\psi\colon X\to X$ such that
$\p\comp\psi$ and $\psi\comp\p$ are chain homotopic to the identity map
$\id_X$.  In this appendix, we show that each complex with injective
components admits a decomposition $X=X'\amalg X''$ such that $X'$ is
homotopically minimal and $X''$ is null homotopic.

Let $\A$ be an abelian category, and suppose that $\A$ admits
injective envelopes. Given a complex $X$ in $\A$ with injective
components, we construct for each $n\in\bbZ$ a new complex $X(n)$ as
follows.  Let $U^n\subseteq X^n$ be the injective envelope of
$Z^nX$. We get a decomposition $X^n=U^n\amalg V^n$. Let $V^{n+1}$ be
the image of $V^{n}$ under the differential $X^n\to X^{n+1}$, and let
$V^p=0$ otherwise. This gives a complex $V$ which is null
homotopic. The canonical map $\iota\colon V\to X$ is a split
monomorphism in each degree. Thus $\iota$ has a left inverse and we
obtain a decomposition $X= U\amalg V$. We put $X(n)=V$.

\begin{lem}\label{le:min}
  Let $\A$ be an abelian category, and suppose that $\A$ admits
  injective envelopes. Then the following are equivalent for a complex
  $X$ in $\A$ with injective components.
\begin{enumerate}
\item The complex $X$ is homotopically minimal.
\item The complex $X$ has no non-zero direct factor which is null homotopic.
\item The canonical map $Z^nX\to X^n$ is an injective envelope for
  all $n\in\bbZ$.
\end{enumerate}
\end{lem}
\begin{proof}
  (1) $\Rightarrow$ (2): Let $X=X'\amalg X''$ and suppose $X'$ is
  null-homotopic. The idempotent map $\e\colon X\to X$ with
  $\Ker\e=X'=\Coker\e$ induces an isomorphism in the homotopy
  category. Thus (1) implies $X'=0$.
  
  (2) $\Rightarrow$ (3): Fix $n\in\bbZ$. Then we have a decomposition
  $X=X(n)\amalg U$ such that $X(n)$ is null homotopic.  Our assumption
  implies $X(n)=0$, and we conclude that the map $Z^nX\to X^n$ is an
  injective envelope.
  
  (3) $\Rightarrow$ (1): Let $\p\colon X\to X$ be a map with inverse
  $\psi$ such that $\psi\comp\p$ and $\p\comp\psi$ are chain homotopic
  to the identity $\id_X$.  Thus we have a family of maps $\r^n\colon
  X^n\to X^{n-1}$ such that
  $$\id_{X^n}=(\psi\comp\p)^n+\d^{n-1}\comp\r^n+\r^{n+1}\comp\d^n.$$
  We claim that $\Ker\p=0$. In fact, we show that
  $K=\Ker(\psi\comp\p)=0$. Let $L^n=K^n\cap Z^nX$. Then $\r^n$
  identifies $L^n$ with $\r^n(L^n)$, and $\r^n(L^n)\cap
  Z^{n-1}X=0$, since $(\d^{n-1}\comp\r^n)L^n=L^n$. The assumption
  on $Z^{n-1}X$ implies $L^n=0$. The same assumption
  on $Z^{n}X$ implies $K^n=0$. Let $C=\Coker\p$. The sequence $0\to
  X\xto{\p}X\to C\to 0$ is split exact in each degree because $X$ has
  injective components. It follows that the sequence is split exact in
  the category of complexes, because $C$ is null homotopic by our
  assumption on $\p$. Let $\p'\colon X\to X$ be a left inverse of
  $\p$. Then $\Ker\p'\cong C$. On the other hand, $\p'$ is
  invertible in the homotopy category of complexes and therefore
  $\Ker\p'=0$ by the first part of this proof. Thus $\p$ is an isomorphism. 
\end{proof}

\begin{prop}\label{pr:min}
  Let $\A$ be an abelian category, and suppose that $\A$ admits
  injective envelopes. Then every complex $X$ in $\A$ with injective
  components has a decomposition $X=X'\amalg X''$ such that $X'$ is
  homotopically minimal and $X''$ is null homotopic. Given a second
  decomposition $X=Y'\amalg Y''$ such that $Y'$ is homotopically
  minimal and $Y''$ is null homotopic, then the canonical map
  $X'\rightarrowtail X\twoheadrightarrow Y'$ is an isomorphism.
\end{prop}
\begin{proof} 
  Take $X''=\coprod_{n\in\bbZ} X(n)$. This complex is null homotopic
  and the canonical map $\iota\colon\coprod_{n\in\bbZ} X(n)\to X$ is a
  split monomorphism in each degree.  Thus $\iota$ has a left inverse
  and we obtain a decomposition $X=X'\amalg X''$. The construction of
  each $X(n)$ shows that the inclusion $Z^n(X')\to (X')^n$ is an injective
  envelope. Thus $X'$ is
  homotopically minimal, by Lemma~\ref{le:min}.

Now let $X=Y'\amalg Y''$ be a second decomposition such that $Y'$ is
  homotopically minimal and $Y''$ is null homotopic. The canonical map
  $\p\colon X'\rightarrowtail X\twoheadrightarrow Y'$ induces an
  isomorphism in the homotopy category, since $X''$ and $Y''$ are null
  homotopic. Thus $\p$ is an isomorphism of complexes, since $X'$ and
  $Y'$ are homotopically minimal. This completes the proof.
\end{proof}

\end{appendix}

\subsection*{Acknowledgements}
It is a pleasure to thank a number of colleagues for their interest in
the topic of this paper: {\O}yvind Solberg helped to prove the crucial
Proposition~\ref{pr:gen}. Peter J{\o}rgensen pointed out the relevance
of the stable derived category in algebraic geometry. Amnon Neeman
communicated an alternative proof of Corollary~\ref{co:products} in
the algebraic geometric context. Ragnar Buchweitz explained the
connection with his earlier unpublished work. Apostolos Beligiannis
provided numerous comments on a preliminary version of this
manuscript. Finally, I wish to thank an anonymous referee for many
helpful suggestions.

\end{document}